\documentclass[final,3p]{elsarticle}
 \usepackage{graphics}
 \usepackage{graphicx}
 \usepackage{epsfig}
\usepackage{amssymb}
 \usepackage{amsthm}
 \usepackage{lineno}
 \usepackage{amsmath}
   \numberwithin{equation}{section}
\usepackage{mathrsfs}

\newtheorem{thm}{Theorem}[section]
\newtheorem{cor}[thm]{Corollary}
\newtheorem{lem}[thm]{Lemma}
\newtheorem{prop}[thm]{Proposition}
\newtheorem{defn}[thm]{Definition}

\biboptions{sort&compress,square}
\begin{document}
\begin{frontmatter}
 \author{Jian Wang}
\author{Yong Wang\corref{cor2}}
\ead{wangy581@nenu.edu.cn}
\cortext[cor2]{Corresponding author}
\address{School of Mathematics and Statistics, Northeast Normal University,
Changchun, 130024, P.R.China}
\title{On the Geometry of Tangent Bundles \\
with the Rescaled Metric}
\begin{abstract}
For a Riemannian manifold $M$, we determine some curvature properties of a tangent bundle equipped with the rescaled metric.
The main aim of this paper is to give explicit formulae for the rescaled metric on $TM$, and investigate the geodesics on the
tangent bundle with respect to the rescaled Sasaki metric.
\end{abstract}
\end{frontmatter}
\section{Introduction}
\label{1}
Tangent boundles of differentiable manifolds are of great importance in many areas of mathematics and physics.
Geometry of the tangent bundle $TM$ of a Riemannian manifold $(M,g)$ with the metric $\bar{g}$ defined by Sasaki in \cite{Sa} had
been studied by many authors. Its construction is based on a natural splitting of the tangent bundle $TTM$ of $TM$ into its
vertical and horizontal subbundles by means of the Levi-Civita connection $\nabla$ on $(M,g)$. The Levi-Civita connection
$\hat{\nabla}$ of the Sasaki metric on $TM$ and its Riemannian curvature tensor $\hat{R}$ were calculated by Kowalski in \cite{Ko}.
With this in hand, the authors derived interesting connections between the geometric properties of $(M,g)$ and $(TM,\hat{g})$ in \cite{Ko}
and \cite{MT}. In \cite{MT}, the authors proved that the Sasaki metric on $TM$ is rather rigid under the scalar curvature
of $(TM,\bar{g})$ is constant.

Another metric nicely fitted to the tangent bundle is the so-called Cheeger-Gromoll metric in \cite{CG}. This can be used to obtain a
natural metric $\tilde{g}$ on the tangent bundle $TM$ of a given Riemannian manifold $(M,g)$. It was expressed more explicitly by
Musso and Tricerri in \cite{MT}. In \cite{Se}, Sekizawa calculated the Levi-Civita connection $\tilde{\nabla}$ and the curvature tensor
$\tilde{R}$ of the tangent bundle $(TM,\tilde{g})$ equipped with the Cheeger-Gromoll metric. Gudmundsson and Kappos derived correct
relations between the geometric properties of $(M,g)$ and $(TM,\tilde{g})$ in \cite{GK1}. In \cite{GK2}, Explicit formulae for the
Cheeger-Gromoll metric on $TM$ was given. The motivation of this paper is to study the geometry of tangent bundles with the
rescaled Sasaki and Cheeger-Gromoll metrics.

This paper is organized as follows: In Section 2, for a Riemannian manifold $(M, g)$, we introduce a natural class of rescaled metrics.
In Section 3, we calculate its Levi-Civita connection, its Riemann curvature tensor associated to the rescaled Sasaki metric.
In Section 4, we investigate geodesics on the tangent bundle with respect to the rescaled Sasaki metric. The main purpose of Section 5
is to obtain some interesting connections between the geometric properties of the manifold $(M, g)$ and its tangent bundle equipped
with the rescaled  Cheeger-Gromoll metric.

\section{Natural Metrics}
In this section we introduce a natural class of rescaled metrics on the tangent bundle $TM$ of a given Riemannian manifold $(M, g)$. This class contains
both the rescaled Sasaki and rescaled Cheeger-Gromoll metrics studied later on.

Throughout this paper we shall assume that $M$ is a smooth $m-$dimensional manifold with maximal atlas
$\mathcal{A}=\{(U_{\alpha},x_{\alpha})|\alpha\in I\}$. For a point $p\in M$, let $T_{p}M$ denote the tangent space of $M$ at $p$.
For local coordinates $(U,x)$ on $M$ and $p\in U$ we define $(\frac{\partial}{\partial x_{k}})_{p}\in T_{p}M$ by
\begin{equation}
(\frac{\partial}{\partial x_{k}})_{p}:f \mapsto \frac{\partial f}{\partial x_{k}}(p)=\partial_{e_{k}}(f\circ x^{-1})(x(p))
\end{equation}
where $\{e_{k}|k=1,\ldots,m\}$ is the standard basis of $\mathbb{R}^{m}$. Then
 $\{(\frac{\partial}{\partial x_{k}})_{p}|k=1,\ldots,m\}$ is a basis for $T_{p}M$. The set $TM=\{(p,u)|p\in M,u\in T_{p}M\}$ is called the
 tangent bundle of $M$ and bundle map $\pi:TM\rightarrow M$ is given by $\pi:(p,u)\mapsto p$.

 As a direct consequence of the Theorem 2.1 in \cite{GK2} we see that the bundle map $\pi:TM\rightarrow M$ is smooth. For each point $p\in M$
 the fiber $\pi^{-1}(p)$ is the tangent space $T_{p}M$ of $M$ at $p$ and hence an $m-$dimensional vector space. For local coordinates
 $(U,x)\in\mathcal{A}$ we define $\bar{x}:\pi^{-1}(U)\rightarrow U\times \mathbb{R}^{m}$ by
\begin{equation}
\bar{x}:(p,\sum_{k=1}^{m}u_{k}\frac{\partial}{\partial x_{k}}|_{p})\mapsto \big(p,(u_{1},\ldots,u_{m})\big).
\end{equation}
The restriction $\bar{x}_{p}=\bar{x}|_{T_{p}M}:T_{p}M\rightarrow \{p\}\times\mathbb{R}^{m}$ to the tangent space $T_{p}M$ is given by
\begin{equation}
\bar{x}_{p}:\sum_{k=1}^{m}u_{k}\frac{\partial}{\partial x_{k}}|_{p}\mapsto (u_{1},\ldots,u_{m})
\end{equation}
so it is obviously a vector space isomorphism. This implies that $\bar{x}:\pi^{-1}(U)\rightarrow U\times \mathbb{R}^{m}$ is
a bundle chart for $TM$. This implies that
\begin{equation}
\mathcal{B}=\{\big(\pi^{-1}(U)\big),\bar{x}|(U,x)\in\mathcal{A}\}
\end{equation}
is a bundle atlas transforming $(TM,M,\pi)$ into an $m-$dimensional topological vector bundle. Since the
manifold $(M,\mathcal{A})$ is smooth the vector bundle $(TM,M,\pi)$ together with the maximal bundle atlas
$\hat{\mathcal{B}}$ induced by $\mathcal{B}$ is a smooth vector bundle.

\begin{defn}
Let $(M, g)$ be a Riemannian manifold. Let $f>0$ and $f\in C^{\infty}(M)$, specially when $f=1$, $\bar{g}^{1}=\bar{g}$.
 A Riemannian rescaled metric $\bar{g}^{f}$ on the
 tangent bundle $TM$ is said to be natural with respect to $g$ on $M$ if
\begin{eqnarray}
 i) \ \bar{g}^{f}_{(p,u)}(X^{h},Y^{h})&=&f(p) g_{p}(X,Y),\\
ii) \ \bar{g}^{f}_{(p,u)}(X^{h}, Y^{v})&=&0
\end{eqnarray}
for all vector fields $X, Y\in C^{\infty}(TM)$ and $(p, u)\in TM.$
\end{defn}
A rescaled natural metric $\bar{g}^{f}$ is constructed in such a way that the vertical and horizontal subbundles are orthogonal and the bundle map
$\pi: (TM, \bar{g}^{f})\rightarrow (M, f g)$ is Riemannian submersion. The rescaled metric $\bar{g}^{f}$ induces a norm on each tangent space of $TM$
which we denote
by $\parallel \cdot \parallel$.
\begin{lem}\label{le:22}
Let $(M, g)$ be a Riemannian manifold and $TM$ be the tangent bundle of $M$. Let $f>0$ and $f\in C^{\infty}(M)$. If the rescaled Riemannian  metric
$\bar{g}^{f}$ on $TM$ is natural with respect to $g$ on $M$ then the corresponding Levi-Civita connection $\overline{\nabla}^{f}$ satisfies
\begin{eqnarray}
i) \  \bar{g}(\overline{\nabla}^{f}_{X^{h}}Y^{h}, Z^{h})&=&\frac{1}{2f}\Big(X(f)g(Y, Z)+Y(f)g(Z, X)-Z(f)g(X, Y)\Big)+g(\nabla_{X}Y,Z),\\
ii) \ \bar{g}(\overline{\nabla}^{f}_{X^{h}}Y^{h}, Z^{v})&=&-\frac{1}{2}\bar{g}\Big((R(X, Y)u)^{v}, Z^{v}\Big) ,\\
iii) \ \bar{g}(\overline{\nabla}^{f}_{X^{h}}Y^{v}, Z^{h})&=&\frac{1}{2f}\bar{g}\Big((R(X, Z)u)^{v}, Y^{v}\Big),\\
iv) \ \bar{g}(\overline{\nabla}^{f}_{X^{h}}Y^{v}, Z^{v})&=&\frac{1}{2}\Big(X^{h}(\bar{g}(Y^{v}, Z^{v}))-\bar{g}(Y^{v}, (\nabla_{X}Z)^{v})
   +\bar{g}(Z^{v}, (\nabla_{X}Y)^{v})\Big),
\end{eqnarray}
\begin{eqnarray}
v) \ \bar{g}(\overline{\nabla}^{f}_{X^{v}}Y^{h}, Z^{h})&=&\frac{1}{2f}\bar{g}\Big((R(Y, Z)u)^{v}, X^{v}\Big),\\
vi) \ \bar{g}(\overline{\nabla}^{f}_{X^{v}}Y^{h}, Z^{v})&=&\frac{1}{2}\Big(Y^{h}(\bar{g}(Z^{v}, X^{v}))-\bar{g}(X^{v}, (\nabla_{Y}Z)^{v})
   -\bar{g}(Z^{v},(\nabla_{Y}X)^{v})\Big),\\
vii) \ \bar{g}(\overline{\nabla}^{f}_{X^{v}}Y^{v}, Z^{h})&=&\frac{1}{2f}\Big(-Z^{h}(\bar{g}(X^{v}, Y^{v}))+\bar{g}(Y^{v}, (\nabla_{Z}X)^{v})
   +\bar{g}(X^{v},(\nabla_{Z}Y)^{v})\Big),\\
viii) \ \bar{g}(\overline{\nabla}^{f}_{X^{v}}Y^{v}, Z^{v})&=&\frac{1}{2}\Big(X^{v}(\bar{g}(X^{v}, Z^{v}))+Y^{v}(\bar{g}(Z^{v}, X^{v}))
   -Y^{v}(\bar{g}(X^{v},Y^{v}))\Big)
\end{eqnarray}
for all vector fields $X, Y, Z\in C^{\infty}(TM)$ and $(p, u)\in TM.$
\end{lem}
\begin{proof}
We shall repeatedly make use of the Kozul formula for the Levi-Civita connection $\overline{\nabla}^{f}$ stating that
\begin{eqnarray}
  2\bar{g}^{f}(\overline{\nabla}^{f}_{X^{i}}Y^{j}, Z^{k})&=&X^{i}(\bar{g}^{f}(Y^{j},Z^{k}))+Y^{j}(\bar{g}^{f}(Z^{k}, X^{i}))
                                                          -Z^{k}(\bar{g}^{f}(X^{i},Y^{j})) \nonumber\\
                                                       &&-\bar{g}^{f}(X^{i},\ [Y^{j}, Z^{k}])+\bar{g}^{f}(Y^{j},[Z^{k},X^{i}])
                                                         +\bar{g}^{f}(Z^{k},\ [X^{i},Y^{j}])
\end{eqnarray}
for all vector fields $X, Y, Z\in\mathcal{C}^{\infty}(TM)$ and $i, j, k\in \{h, v\}.$

$i)$ The result is a direct consequence of the following calculations using Definition 2.1 and Proposition 5.1 in \cite{GK2},
\begin{eqnarray}
  2\bar{g}^{f}(\overline{\nabla}^{f}_{X^{h}}Y^{h}, Z^{h})&=&X^{h}(\bar{g}^{f}(Y^{h}, Z^{h}))+Y^{h}(\bar{g}^{f}(Z^{h}, X^{h}))
                                                          -Z^{h}(\bar{g}^{f}(X^{h}, Y^{h}))\nonumber\\
                                                       &&-\bar{g}^{f}(X^{h}, [Y^{h}, Z^{h}])+\bar{g}^{f}(Y^{h}, [Z^{h}, X^{h}])
                                                         +\bar{g}^{f}(Z^{h}, [X^{h}, Y^{h}])\nonumber\\
                                                       &=&X^{h}(f g(Y,Z)\circ\pi)+Y^{h}(f g(Z, X)\circ\pi)
                                                          -Z^{h}(f g(X, Y)\circ\pi)\nonumber\\
                                                       &&-\bar{g}^{f}(X^{h}, [Y, Z]^{h})+\bar{g}^{f}(Y^{h}, [Z, X]^{h})
                                                         +\bar{g}^{f}(Z^{h},[X, Y]^{h})\nonumber\\
                                                       &=& X(f)g(Y, Z)+Y(f)g(Z, X)-Z(f)g(X,  Y)+2f\bar{g}^{f}(\nabla_{X}Y),  Z).
\end{eqnarray}

$ii)$ The statement is obtained as follows.
\begin{eqnarray}
  2\bar{g}^{f}(\overline{\nabla}^{f}_{X^{h}}Y^{v}, Z^{h})&=&X^{h}(\bar{g}^{f}(Y^{h},Z^{v}))+Y^{h}(\bar{g}^{f}(Z^{v}, X^{h}))
                                                          -Z^{v}(\bar{g}^{f}(X^{h}, Y^{h})) \nonumber\\
                                                       &&-\bar{g}^{f}(X^{h},[Y^{h}, Z^{v}])+\bar{g}^{f}(Y^{h}, [Z^{v}, X^{h}])
                                                         +\bar{g}^{f}(Z^{v}, [X^{h}, Y^{h}]) \nonumber\\
                                                       &=&-Z^{v}(f g(X, Y))+\bar{g}^{f}(Z^{v}, [X^{h},Y^{h}])  \nonumber\\
                                                       &=&-\bar{g}^{f}(Z^{v},(R(X, Y)u)^{v})
\end{eqnarray}

$iii)$ and $v)$ are analogous to $ii)$.

$iv)$ Again using Definition 2.1 and Proposition 5.1 in \cite{GK2} we yield
\begin{eqnarray}
  2\bar{g}^{f}(\overline{\nabla}^{f}_{X^{h}}Y^{v}, Z^{v})&=&X^{h}(\bar{g}^{f}(Y^{v}, Z^{v}))+Y^{v}(\bar{g}^{f}(Z^{v}, X^{h}))
                                                          -Z^{v}(\bar{g}^{f}(X^{h}, Y^{v})) \nonumber\\
                                                       &&-\bar{g}^{f}(X^{h}, [Y^{v}, Z^{v}])+\bar{g}^{f}(Y^{v}, [Z^{v}, X^{h}])
                                                         +\bar{g}^{f}(Z^{v}, [X^{h}, Y^{v}]) \nonumber\\
                                                       &=&X^{h}(\bar{g}(Y^{v}, Z^{v}))-\bar{g}(Y^{v},(\nabla_{X}Z)^{v})
                                                         +\bar{g}(Z^{v}, (\nabla_{X}Y)^{v})
\end{eqnarray}

$vi)$ and $vii)$ are analogous to iv).

$viii)$ The statement is a direct consequence of the fact that the Lie bracket of two vertical vector fields vanishes.
\end{proof}
\begin{cor}\label{co:23}
Let $(M, g)$ be a Riemannian manifold and $\bar{g}^{f}$ be a rescaled natural rescaled metric on the tangent bundle $TM$ of $M$.
 Then the Levi-Civita connection
$\overline{\nabla}^{f}$ satisfies
\begin{equation}
(\overline{\nabla}^{f}_{X^{h}}Y^{h})_{(p, u)}=(\nabla^{f}_{X}Y)^{h}_{(p, u)}-\frac{1}{2}\Big(R(X, Y)u\Big)^{v}+
                                              \frac{1}{2f(p)}\Big(X(f)Y+Y(f)X-g(X, Y)\circ\pi(\texttt{d}(f\circ\pi))^{*}\Big)^{h}_{p}
\end{equation}
for all vector fields $X, Y\in C^{\infty}(TM)$ and $(p, u)\in TM.$
\end{cor}
\begin{proof}
By proposition 3.5 in \cite{GK2}, each tangent vector $Z\in T_{(p,u)}TM$ can be decomposed as $Z=Z^{h}_{1}+Z^{v}_{2}$.
Using $i)$ and $ii)$ of Lemma 2.2, we have
\begin{eqnarray}
  \bar{g}(\overline{\nabla}^{f}_{X^{h}}Y^{h}, Z^{h}_{1}+Z^{v}_{2})&=&-\frac{1}{2}\bar{g}\Big((R(X, Y)u)^{v}, Z^{h}_{1}+Z^{v}_{2}\Big)
                                                 +g\Big((\nabla_{X}Y)^{h},Z^{h}_{1}+Z^{v}_{2}\Big)\nonumber\\
                                              &&+\frac{1}{2f}\Big(X(f)g(Y^{h}, Z^{h}_{1}+Z^{v}_{2})+g((Y f X)^{h}, Z^{h}_{1}+Z^{v}_{2})\nonumber\\
                                              && -\langle g(X^{h}, Y^{h})\texttt{d}(f\circ\pi), Z^{h}_{1}+Z^{v}_{2}\rangle\Big) \nonumber\\
                                              &=&(\nabla^{f}_{X}Y)^{h}-\frac{1}{2}\Big(R(X, Y)u\Big)^{v}+\frac{1}{2f}\Big(X(f)Y+Y(f)X\nonumber\\
                                             &&-g(X, Y)\circ\pi(\texttt{d}(f\circ\pi))^{*}\Big)^{h}.
\end{eqnarray}
\end{proof}

\begin{defn}
Let $(M, g)$ be a Riemannian manifold and $F:TM\rightarrow TM$ be a smooth bundle endomorphism of the tangent bundle $TM$. Then we define
the vertical and horizontal lifts $F^{v}:TM\rightarrow TTM$, $F^{h}:TM\rightarrow TTM$ of $F$ by
\begin{equation}
F^{v}(\eta)=\sum_{i=1}^{m}\eta_{i}F(\partial_{i})^{v}  \quad and \quad F^{h}(\eta)=\sum_{i=1}^{m}\eta_{i}F(\partial_{i})^{h},
\end{equation}
where $\sum_{i=1}^{m}\eta_{i}\partial_{i}\in\pi^{-1}(V)$ is a local representation of $\eta\in C^{\infty}(TM)$.
\end{defn}
\begin{lem}\label{le:25}
Let $(M, g)$ be a Riemannian manifold and the tangent bundle $TM$ be equipped with a rescaled metric $\bar{g}^{f}$ which is natural with respect to $g$
on $M$. If $F:TM\rightarrow TM$ is a smooth bundle endomorphism of the tangent bundle, then
\begin{eqnarray}
i) \ (\overline{\nabla}^{f}_{X^{v}}F^{v})_{\xi}&=&F(X_{p})^{v}_{\xi}+\sum_{i=1}^{m}u(x_{i})(\overline{\nabla}^{f}_{X^{v}}F(\partial_{i})^{v})_{\xi},\\
ii) \ (\overline{\nabla}^{f}_{X^{v}}F^{h})_{\xi}&=&F(X_{p})^{h}_{\xi}+\sum_{i=1}^{m}u(x_{i})(\overline{\nabla}^{f}_{X^{v}}F(\partial_{i})^{h})_{\xi},\\
iii) \ (\overline{\nabla}^{f}_{X^{h}}F^{v})_{\xi}&=&(\overline{\nabla}^{f}_{X^{h}}F(u)^{v})_{\xi},\\
iv) \ (\overline{\nabla}^{f}_{X^{h}}F^{h})_{\xi}&=&(\overline{\nabla}^{f}_{X^{h}}F(u)^{h})_{\xi},
\end{eqnarray}
for any $X\in C^{\infty}(TM)$, $\xi=(p, u)\in TM$ and $\eta=\sum_{i=1}^{m}\eta_{i}\partial_{i}\in\pi^{-1}(V)$.
\end{lem}
\begin{proof}
Let $(x_{1},\cdots,x_{m})$ be local coordinates on $M$ in a neighborhood $V$ of $p$. Then, using the abbreviation $X_{i}$ for
$\frac{\partial}{\partial x_{i}}$, we have $X^{v}(\texttt{d}x_{i})=\texttt{d}x_{i}(X)=X(x_{i})$ and $\texttt{d}x_{i}(p,u)=\eta_{i}(p)$
for $i\in\{1,\cdots,m\}$. Hence
\begin{eqnarray}
  (\overline{\nabla}^{f}_{X^{v}}F^{v})_{\xi}&=&\sum_{i=1}^{m}\overline{\nabla}^{f}_{X^{v}}(\eta_{i}F(\partial_{i})^{v})
                                            =\sum_{i=1}^{m}X^{v}(\texttt{d}x_{i})F(\partial_{i})^{v}
                                              +\eta_{i}\overline{\nabla}^{f}_{X^{v}}F(\partial_{i})^{v}  \nonumber\\
                                            &=&\sum_{i=1}^{m}X(x_{i})F(\partial_{i})^{v}
                                              +\eta_{i}\overline{\nabla}^{f}_{X^{v}}F(\partial_{i})^{v}
                                              =F(X_{p})^{v}_{\xi}+\sum_{i=1}^{m}u(x_{i})(\overline{\nabla}^{f}_{X^{v}}F(\partial_{i})^{v})_{\xi}.
\end{eqnarray}
Similarly we have
\begin{eqnarray}
  (\overline{\nabla}^{f}_{X^{v}}F^{h})_{\xi}&=&\sum_{i=1}^{m}\overline{\nabla}^{f}_{X^{v}}(\eta_{i}F(\partial_{i})^{h})
                                            =\sum_{i=1}^{m}X^{v}(\texttt{d}x_{i})F(\partial_{i})^{h}
                                              +\eta_{i}\overline{\nabla}^{f}_{X^{v}}F(\partial_{i})^{h} \nonumber\\
                                            &=&\sum_{i=1}^{m}X(x_{i})F(\partial_{i})^{h}
                                              +\eta_{i}\overline{\nabla}^{f}_{X^{v}}F(\partial_{i})^{h}
                                              =F(X_{p})^{h}_{\xi}+\sum_{i=1}^{m}u(x_{i})(\overline{\nabla}^{f}_{X^{v}}F(\partial_{i})^{h})_{\xi}.
\end{eqnarray}
For the last two equations of the lemma we use a differentiable curve $\gamma: [0, 1]\rightarrow M$ such that $\gamma(0)=p$ and $\gamma'(0)=X_{p}$
to get a differentiable curve $U\circ\gamma: [0, 1]\rightarrow TM$ such that $U\circ\gamma(0)=\xi$ and $(U\circ\gamma)'(0)=X^{h}_{\xi}$. By the
definition of $F^{v}$ and $F^{h}$ we get
\begin{eqnarray}
  F^{v}|_{U\circ\gamma(t)}&=&\sum_{i=1}^{m}\texttt{d}x_{i}F(\partial_{i})^{v}|_{U\circ\gamma(t)}
                       =\sum_{i=1}^{m}\texttt{d}x_{i}(U\circ\gamma(t))F(\partial_{i})^{v}|_{U\circ\gamma(t)}   \nonumber\\
                       &=&F(\sum_{i=1}^{m}u(x_{i})_{p}e_{j})^{v}|_{U\circ\gamma(t)}=(F\circ U)^{v}|_{U\circ\gamma(t)}.
\end{eqnarray}
Similarly $F^{h}|_{U\circ\gamma}=(F\circ U)^{h}_{U\circ\gamma}$. This proves parts $iii)$ and $iv)$.
\end{proof}

\section{The Rescaled Sasaki Metric}
This section is devoted to the Sasaki metric $\hat{g}$ on the tangent bundle $TM$ introduced by Sasaki in the famous paper \cite{Sa}.
We calculate its Levi-Civita connection $\hat{\nabla}^{f}$, its Riemann curvature tensor and obtain some interesting connections between the
geometric properties of the manifold $(M, g)$ and its tangent bundle $(TM, \hat{g}^{f})$ equipped with the rescaled Sasaki metric.
\begin{defn}
Let $(M, g)$ be a Riemannian manifold. Let $f>0$ and $f\in C^{\infty}(M)$. Then the rescaled Sasaki metric $\hat{g}^{f}$ on the tangent
bundle $TM$ of $M$ is given by
\begin{eqnarray}
i) \ \hat{g}_{(x,u)}^{f}(X^{h}, Y^{h})&=&f(p)g_{p}(X, Y),\\
ii) \ \hat{g}_{(x,u)}^{f}(X^{v},  Y^{h})&=&0 ,\\
iii) \ \hat{g}_{(x,u)}^{f}(X^{v},  Y^{v})&=&g_{p}(X, Y).
\end{eqnarray}
for all vector fields $X, Y\in C^{\infty}(TM)$.
\end{defn}
The rescaled Sasaki metric is obviously contained in the class of rescaled $g-$natural metrics. It is constructed in such a manner that inner
products are respected not only by lifting vectors horizontally but vertically as well.

\begin{prop}\label{pr: 32}
Let $(M, g)$ be a Riemannian manifold and $\hat{\nabla}^{f}$ be Levi-Civita connection of the tangent bundle $(TM, \hat{g}^{f})$ equipped with
the rescaled Sasaki metric. Then
\begin{eqnarray}
i) \ (\hat{\nabla}^{f}_{X^{h}}Y^{h})_{(p,u)}&=&(\nabla_{X}Y)^{h}_{(p,u)}+\frac{1}{2f(p)}
          \Big((X(f)Y+Y(f)X)-g(X,Y)\circ\pi(\texttt{d}(f\circ\pi))^{*}\Big)^{h}_{p}\nonumber\\
&&-\frac{1}{2}\Big(R_{p}(X,Y)u\Big)^{v},\\
ii) \ (\hat{\nabla}^{f}_{X^{h}}Y^{v})_{(p,u)}&=&(\nabla_{X}Y)^{v}_{(p,u)}+\frac{1}{2f(p)}\Big(R_{p}(u,Y)X\Big)^{h},\\
iii) \ (\hat{\nabla}^{f}_{X^{v}}Y^{h})_{(p,u)}&=&\frac{1}{2f(p)}\Big(R_{p}(u,X)Y\Big)^{h},\\
iv) \ (\hat{\nabla}^{f}_{X^{v}}Y^{v})_{(p,u)}&=&0
\end{eqnarray}
for any $X,Y\in C^{\infty}(TM)$, $\xi=(p, u)\in TM$.
\end{prop}
\begin{proof}
 $i)$ The statement is a direct consequence of Corollary 2.3.

 $ii)$ By applying Lemma 2.2 we obtain the following for the horizontal part
\begin{eqnarray}
  2\hat{g}^{f}(\hat{\nabla}^{f}_{X^{h}}Y^{v}, Z^{h})&=&-\hat{g}^{f}((R(Z,X)u)^{v},Y^{v})=-g(R(u,Y)Z,X)\nonumber\\
                                              &=&g(R(u,Y)X,Z)=\frac{1}{f}\hat{g}^{f}((R(u,Y)X)^{h},Z^{h}),
\end{eqnarray}
As for the vertical part note that
\begin{eqnarray}
  2\hat{g}^{f}(\hat{\nabla}^{f}_{X^{h}}Y^{v}, Z^{v})&=&X^{h}(\hat{g}^{f}(Y^{v},Z^{v}))+\hat{g}^{f}(Z^{v},(\nabla_{X}Y)^{v})
                                                          -\hat{g}^{f}(Y^{v},(\nabla_{X}Z)^{v})  \nonumber\\
                                              &=&X(g(Y,Z))+g(Z,\nabla_{X}Y)-g(Y,\nabla_{X}Z)     \nonumber\\
                                              &=&2\hat{g}^{f}((\nabla_{X}Y)^{v},Z^{v}).
\end{eqnarray}

$iii)$ For the horizontal part we get calculations similar to those above
\begin{eqnarray}
  2\hat{g}(\hat{\nabla}^{f}_{X^{v}}Y^{h}, Z^{h})&=&\frac{1}{f}\hat{g}(X^{v},(R(Y,Z)u)^{v})=\frac{1}{f}g(X,R(Y,Z)u)\nonumber\\
                                              &=&\frac{1}{f}g(R(u,X)Y,Z).
\end{eqnarray}
The rest follows by
\begin{eqnarray}
  2\hat{g}(\hat{\nabla}^{f}_{X^{v}}Y^{h}, Z^{v})&=&Y^{h}(\hat{g}(Z^{v},X^{v}))-\hat{g}(Z^{v},(\nabla_{Y}X)^{v})
                                                          -\hat{g}(X^{v},(\nabla_{Y}Z)^{v})  \nonumber\\
                                              &=&Y(g(Z,X))-g(Z,\nabla_{Y}X)-g(X,\nabla_{Y}Z)=0.
\end{eqnarray}

$iv)$ Using Lemma 2.2 again we yield
\begin{eqnarray}
  2f\hat{g}(\hat{\nabla}^{f}_{X^{v}}Y^{v}, Z^{h})&=&-Z^{h}(\hat{g}(X^{v},Y^{v}))+\hat{g}(Y^{v},(\nabla_{Z}X)^{v})
                                                         +\hat{g}(X^{v},(\nabla_{Z}Y)^{v}) \nonumber\\
                                              &=&-Z(g(X,Y))+g(Y,\nabla_{Z}X)+g(X,\nabla_{Z}Y)=0,
\end{eqnarray}
and
\begin{eqnarray}
  2\hat{g}(\hat{\nabla}^{f}_{X^{v}}Y^{v}, Z^{v})&=&X^{v}(\hat{g}(Y^{v},Z^{v}))+Y^{v}(\hat{g}(Z^{v},X^{v}))
                                                     -Z^{v}(\hat{g}(X^{v},Y^{v}))  \nonumber\\
                                              &=&X^{v}(g(Y,Z))+Y^{v}(g(Z,X))-Z^{v}(g(X,Y))=0.
\end{eqnarray}
This completes the proof.
\end{proof}

We shall now turn our attention to the Riemann Curvature tensor $\hat{R}^{f}$ of the tangent bundle $TM$ equipped with the
rescaled Sasaki metric $\hat{g}^{f}$. For this we need the following useful Lemma.
\begin{lem}\label{le:33}
Let $(M, g)$ be a Riemannian manifold and $\hat{\nabla}^{f}$ be the Levi-Civita connection of the tangent bundle $(TM, \hat{g}^{f})$,
equipped with the rescaled Sasaki metric $\hat{g}^{f}$. Let $F:TM\rightarrow TM$ is a smooth bundle endomorphism of the tangent bundle, then
\begin{equation}
(\hat{\nabla}^{f}_{X^{v}}F^{v})_{\xi}=F(X_{p})^{v}_{\xi},
\end{equation}
and
\begin{equation}
(\hat{\nabla}^{f}_{X^{v}}F^{h})_{\xi}=F(X_{p})^{h}_{\xi}+\frac{1}{2f(p)}\Big(R(u,X)F(u)\Big)^{h}_{\xi}
\end{equation}
for any $X\in C^{\infty}(TM)$ and $\xi=(p, u)\in TM$.
\end{lem}
\begin{proof}
By applying $i)$ of Lemma 2.5 and $iv)$ of Proposition 3.2 we obtain the following
\begin{equation}
(\hat{\nabla}^{f}_{X^{v}}F^{v})_{\xi}=F(X_{p})^{v}_{\xi}+\sum_{i=1}^{m}u(x_{i})(\overline{\nabla}^{f}_{X^{v}}F(\partial_{i})^{v})_{\xi}
                                          =F(X_{p})^{v}_{\xi}.
\end{equation}
By applying $ii)$ of Lemma 2.5 and $iii)$ of Proposition 3.2, we get
\begin{equation}
(\hat{\nabla}^{f}_{X^{v}}F^{h})_{\xi}=F(X_{p})^{h}_{\xi}+(\overline{\nabla}^{f}_{X^{v}}F(u)^{h})_{\xi}
                                          =F(X_{p})^{h}_{\xi}+\frac{1}{2f(p)}\Big(R(u,X)F(u)\Big)^{h}_{\xi}.
\end{equation}

\end{proof}
\begin{prop}\label{pr: 34}
Let $(M, g)$ be a Riemannian manifold and $\hat{R}^{f}$ be the Riemann curvature tensor of the tangent bundle $(TM, \hat{g}^{f})$ equipped with
the rescaled Sasaki metric. Then the following formulae hold
\begin{eqnarray}
i) \ \hat{R}^{f}_{(p,u)}(X^{v},Y^{v})Z^{v}&=&0,\\
ii)\ \hat{R}^{f}_{(p,u)}(X^{h},Y^{v})Z^{v}&=&\Big(-\frac{1}{2f(p)}R(Y,Z)X-\frac{1}{4f^{2}(p)}R(u,Y)(R(u,Z)X)\Big)^{h}_{p},
\end{eqnarray}
\begin{eqnarray}
iii)\ \hat{R}^{f}_{(p,u)}(X^{v},Y^{v})Z^{h}&=&\Big(-\frac{1}{2f(p)}R(Y,X)Z-\frac{1}{4f^{2}(p)}R(u,Y)(R(u,X)Z)\Big)^{h}_{p} \nonumber\\
                                           &&+\Big(\frac{1}{2f(p)}R(X,Y)Z+\frac{1}{4f^{2}(p)}R(u,y)(R(u,Y)Z)\Big)^{h}_{p},
\end{eqnarray}

\begin{eqnarray}
iv) \ \hat{R}^{f}_{(p,u)}(X^{h},Y^{v})Z^{h}&=&\Big(\nabla_{X}(\frac{1}{2f(p)}R(u,Y)Z)\Big)^{h}_{p}+\frac{1}{4f(p)}R((R(u,Y)Z),X)u\nonumber\\
                  &&+A_{f}\Big(X,\frac{1}{2f(p)}(R(u,Y)Z)\Big) \nonumber\\
                  &&-\frac{1}{2f(p)}\Big(R(u,Y)(\nabla_{X}Z+A_{f}(X,Z))\Big)^{h}_{p} \nonumber\\
                  &&+\frac{1}{2}\Big(R(X,Z)u\Big)^{v}_{p}-\frac{1}{2f(p)}\Big(R(u,\nabla_{X}Y)Z\Big)^{h}_{p},
\end{eqnarray}
\begin{eqnarray}
v) \ \hat{R}^{f}_{(p,u)}(X^{h},Y^{h})Z^{v}&=&\Big(\nabla_{X}(\frac{1}{2f(p)}R(u,Z)Y)\Big)^{h}_{p}-\Big(\nabla_{Y}
                                              (\frac{1}{2f(p)}R(u,Z)X)\Big)^{h}_{p}\nonumber\\
                                             && +\frac{1}{4f(p)}R(R(u,Z)Y,X)u-\frac{1}{4f(p)}R(R(u,Z)X,Y)u \nonumber\\
                                             &&+\frac{1}{2f(p)}A_{f}(X,R(u,Z)Y)-\frac{1}{2f(p)}A_{f}(Y,R(u,Z)X) \nonumber\\
                                             &&+\frac{1}{2f(p)}R(u,Z)[Y,X]+\Big(R(X,Y)u\Big)^{v}_{p}\nonumber\\
                                             &&+\frac{1}{2f(p)}\Big(R(u,\nabla_{Y}Z)X\Big)^{h}_{p}-\frac{1}{2f(p)}\Big(R(u,\nabla_{X}Z)Y\Big)^{h}_{p},
\end{eqnarray}
\begin{eqnarray}
vi) \ \hat{R}^{f}_{(p,u)}(X^{h},Y^{h})Z^{h} &=&\hat{\nabla}^{f}_{X^{h}}\hat{\nabla}^{f}_{Y^{h}}Z^{h}-\hat{\nabla}^{f}_{Y^{h}}\hat{\nabla}^{f}_{X^{h}}Z^{h}
                                 -\hat{\nabla}^{f}_{[X^{h},Y^{h}]}Z^{h}\nonumber\\
                             &=&\hat{\nabla}^{f}_{X^{h}}(F_{1}^{h})-\hat{\nabla}^{f}_{Y^{h}}\Big((\nabla_{X}Z)^{h}
                                +A_{f}(X,Z)^{h}+F_{2}^{h}\Big)-\hat{\nabla}^{f}_{(\nabla_{X}Y)^{h}}Z^{h} \nonumber\\
                             &=&\nabla_{X}\Big(\nabla_{Y}Z+A_{f}(Y,Z)\Big)^{h}+A_{f}\Big(X,\nabla_{Y}Z+A_{f}(Y,Z)\Big)^{h}\nonumber\\
                                           &&-\frac{1}{2}\Big(R(X,\nabla_{Y}Z+A_{f}(Y,Z))u\Big)^{v}
                                              -\nabla_{Y}\Big(\nabla_{X}Z+A_{f}(X,Z)\Big)^{h} \nonumber\\
                                           &&-A_{f}\Big(Y,\nabla_{X}Z+A_{f}(X,Z)\Big)^{h}+\frac{1}{2}\Big(R(Y,\nabla_{X}Z+A_{f}(X,Z))u\Big)^{v} \nonumber\\
                                           &&-\Big(\nabla_{[X,Y]}Z\Big)^{h}-A_{f}([X,Y],Z)^{h}-\frac{1}{2}\Big(R([X,Y],Z)u\Big)^{v}\nonumber\\
                                           &&+\frac{1}{2f}\Big(R(u,R(X,Y)u)Z\Big)^{h}+\frac{1}{2}\Big(\nabla_{Y}(R(X,Z)u)\Big)^{v}\nonumber\\
                                           &&+\frac{1}{4f}\Big(R(u,R(X,Z)u)Y\Big)^{h}-\frac{1}{2}\Big(\nabla_{X}(R(Y,Z)u)\Big)^{v}\nonumber\\
                                           &&-\frac{1}{4f}\Big(R(u,R(Y,Z)u)X\Big)^{h}.
\end{eqnarray}
for any $X, Y, Z\in T_{p}M$.
\end{prop}
\begin{proof}
$i)$ The result follows immediately from Proposition 3.2.

$ii)$ Let $F: TM\rightarrow TM$ be the bundle endomorphism given by
\begin{equation}
F: u\mapsto \frac{1}{2f}R(u,Z)X.
\end{equation}
 Applying Proposition 3.2 and Lemma 3.3 we have
\begin{equation}
\hat{\nabla}^{f}_{Y^{v}}F^{h}=F(Y)^{h}+\frac{1}{2f}\Big(R(u,Y)F(u)\Big)^{h}.
\end{equation}
This implies that

\begin{eqnarray}
\hat{R}^{f}(X^{h},Y^{v})Z^{v}&=&\hat{\nabla}^{f}_{X^{h}}\hat{\nabla}^{f}_{Y^{v}}Z^{v}-\hat{\nabla}^{f}_{Y^{v}}\hat{\nabla}^{f}_{X^{h}}Z^{v}
                                           -\hat{\nabla}^{f}_{[X^{h},Y^{v}]}Z^{v}\nonumber\\
                                   &=&-\hat{\nabla}^{f}_{Y^{v}}\hat{\nabla}^{f}_{X^{h}}Z^{v}
                                        =-\hat{\nabla}^{f}_{Y^{v}}\Big((\nabla_{X}Z)^{v}+F^{h}\Big)\nonumber\\
                                   &=&-\hat{\nabla}^{f}_{Y^{v}}F^{h}=-F(Y)^{h}-\frac{1}{2f}\Big(R(u,Y)F(u)\Big)^{h} \nonumber\\
                                   &=&\Big(-\frac{1}{2f}R(Y,Z)X-\frac{1}{4f^{2}}R(u,Y)(R(u,Z)X)\Big)^{h}.
\end{eqnarray}

$iii)$ Using $ii)$ and $1^{st}$ Bianchi identity we get
\begin{equation}
\hat{R}^{f}(X^{v},Y^{v})Z^{h}=\hat{R}^{f}(Z^{h},Y^{v})X^{v}-\hat{R}^{f}(Z^{h},X^{v})Y^{v}
\end{equation}
which gives
\begin{eqnarray}
 \hat{R}^{f}(X^{v},Y^{v})Z^{h}&=&\Big(-\frac{1}{2f}R(Y,X)Z-\frac{1}{4f^{2}}R(u,Y)(R(u,X)Z)\Big)^{h} \nonumber\\
                                           &&+\Big(\frac{1}{2f}R(X,Y)Z+\frac{1}{4f^{2}}R(u,X)(R(u,Y)Z)\Big)^{h}.
\end{eqnarray}

$iv)$ Let $F_{1}, F_{2}: TM\rightarrow TM$ be the bundle endomorphisms given by
\begin{equation}
F_{1}(u)\mapsto \frac{1}{2f}R(u,Y)Z \quad and \quad F_{2}(u)\mapsto -\frac{1}{2f}R(X,Z)u.
\end{equation}
Then Proposition 3.2 implies that
\begin{eqnarray}
\hat{R}^{f}(X^{h},Y^{v})Z^{h}&=&\hat{\nabla}^{f}_{X^{h}}\hat{\nabla}^{f}_{Y^{v}}Z^{h}-\hat{\nabla}^{f}_{Y^{v}}\hat{\nabla}^{f}_{X^{h}}Z^{h}
                                 -\hat{\nabla}^{f}_{[X^{h},Y^{v}]}Z^{h}\nonumber\\
                             &=&\hat{\nabla}^{f}_{X^{h}}(F_{1}^{h})-\hat{\nabla}^{f}_{Y^{v}}\Big((\nabla_{X}Z)^{h}
                                +A_{f}(X,Z)^{h}+F_{2}^{v}\Big)-\hat{\nabla}^{f}_{(\nabla_{X}Y)^{v}}Z^{h}\nonumber\\
                             &=&(\nabla_{X}F_{1}(u))^{h}-\frac{1}{2}\Big(R(X,F_{1}(u))u\Big)^{v}+A_{f}(X,F_{1}(u))^{h}\nonumber\\
                                &&-\frac{1}{2f}\Big(R(u,Y)(\nabla_{X}Z+A_{f}(X,Z))\Big)^{h}-F_{2}(Y)^{v}
                                -\frac{1}{2f}\Big(R(u,\nabla_{X}Y)Z\Big)^{h}\nonumber\\
                            &=&\Big(\nabla_{X}(\frac{1}{2f}R(u,Y)Z)\Big)^{h}+\frac{1}{4f}R(R(u,Y)Z,X)u\nonumber\\
                              &&+A_{f}\Big(X,\frac{1}{2f}(R(u,Y)Z)\Big)-\frac{1}{2f}\Big(R(u,Y)(\nabla_{X}Z+A_{f}(X,Z))\Big)^{h} \nonumber\\
                              &&+\frac{1}{2}\Big(R(X,Z)u\Big)^{v}-\frac{1}{2f}\Big(R(u,\nabla_{X}Y)Z\Big)^{h}.
\end{eqnarray}

$v)$ Applying part $iv)$ and $1^{st}$ Bianchi identity
\begin{equation}
\hat{R}^{f}(X^{h},Y^{h})Z^{v}=\hat{R}^{f}(X^{h},Z^{v})Y^{h}-\hat{R}^{f}(Y^{h},Z^{v})X^{h},
\end{equation}
we get
\begin{eqnarray}
\hat{R}^{f}(X^{h},Y^{h})Z^{v}&=&\Big(\nabla_{X}(\frac{1}{2f}R(u,Z)Y)\Big)^{h}+\frac{1}{4f}R(R(u,Z)Y,X)u+A_{f}\Big(X,\frac{1}{2f}(R(u,Z)Y)\Big)\nonumber\\
                             &&-\frac{1}{2f}\Big(R(u,Z)(\nabla_{X}Y+A_{f}(X,Y))\Big)^{h}+\frac{1}{2}(R(X,Y)u)^{v}
                                -\frac{1}{2f}\Big(R(u,\nabla_{X}Z)Y\Big)^{h}\nonumber\\
                             &&-\Big(\nabla_{Y}(\frac{1}{2f}R(u,Z)X)\Big)^{h}-\frac{1}{4f}R(R(u,Z)X,Y)u-A_{f}\Big(Y,\frac{1}{2f}R(u,Z)X\Big)\nonumber\\
                             &&+\frac{1}{2f}\Big(R(u,Z)(\nabla_{Y}X+A_{f}(Y,X))\Big)^{h}-\frac{1}{2}(R(Y,X)u)^{v}
                                +\frac{1}{2f}\Big(R(u,\nabla_{Y}Z)X\Big)^{h},\nonumber\\
\end{eqnarray}
from which the result follows.

$vi)$ By $i)$ of Proposition 3.2 and direct calculation we get
\begin{eqnarray}
\hat{R}^{f}(X^{h},Y^{h})Z^{h}&=&\hat{\nabla}^{f}_{X^{h}}\hat{\nabla}^{f}_{Y^{h}}Z^{h}-\hat{\nabla}^{f}_{Y^{h}}\hat{\nabla}^{f}_{X^{h}}Z^{h}
                                 -\hat{\nabla}^{f}_{[X^{h},Y^{h}]}Z^{h}\nonumber\\
                             &=&\hat{\nabla}^{f}_{X^{h}}(F_{1}^{h})-\hat{\nabla}^{f}_{Y^{h}}\Big((\nabla_{X}Z)^{h}
                                +A_{f}(X,Z)^{h}+F_{2}^{h}\Big)-\hat{\nabla}^{f}_{(\nabla_{X}Y)^{h}}Z^{h} \nonumber\\
                             &=&\nabla_{X}\Big(\nabla_{Y}Z+A_{f}(Y,Z)\Big)^{h}+A_{f}\Big(X,\nabla_{Y}Z+A_{f}(Y,Z)\Big)^{h}\nonumber\\
                                           &&-\frac{1}{2}\Big(R(X,\nabla_{Y}Z+A_{f}(Y,Z))u\Big)^{v}
                                              -\nabla_{Y}\Big(\nabla_{X}Z+A_{f}(X,Z)\Big)^{h} \nonumber\\
                                           &&-A_{f}\Big(Y,\nabla_{X}Z+A_{f}(X,Z)\Big)^{h}+\frac{1}{2}\Big(R(Y,\nabla_{X}Z+A_{f}(X,Z))u\Big)^{v} \nonumber\\
                                           &&-\Big(\nabla_{[X,Y]}Z\Big)^{h}-A_{f}([X,Y],Z)^{h}-\frac{1}{2}\Big(R([X,Y],Z)u\Big)^{v}\nonumber\\
                                           &&+\frac{1}{2f}\Big(R(u,R(X,Y)u)Z\Big)^{h}+\frac{1}{2}\Big(\nabla_{Y}(R(X,Z)u)\Big)^{v}\nonumber\\
                                           &&+\frac{1}{4f}\Big(R(u,R(X,Z)u)Y\Big)^{h}-\frac{1}{2}\Big(\nabla_{X}(R(Y,Z)u)\Big)^{v} \nonumber\\
                                           &&-\frac{1}{4f}\Big(R(u,R(Y,Z)u)X\Big)^{h}.
\end{eqnarray}
\end{proof}

We shall now compare the geometries of the manifold $(M,g)$ and its tangent bundle $TM$ equipped with the rescaled Sasaki metric $\hat{g}^{f}$.
\begin{thm}\label{th:35}
Let $(M, g)$ be a Riemannian manifold and $TM$ be its tangent bundle with the rescaled Sasaki metric $\hat{g}^{f}$. Then $TM$ is flat if and only if
$M$ is flat and $f=C(constant)$.
\end{thm}
\begin{proof}
Applying  proposition 3.4 and
\begin{equation}
A_{f}(X,Y)=\frac{1}{2f}\Big(X(f)Y+Y(f)X-g(X,Y)(\texttt{d}f)^{*}\Big)^{h}.
\end{equation}
If $A_{f}=0$, we have
\begin{equation}
X(f)Y+Y(f)X-g(X,Y)(\texttt{d}f)^{*}=0,
\end{equation}
then $R\equiv 0$ implies $\hat{R}^{f}\equiv 0$. If we assume that $\hat{R}^{f}\equiv 0$ and calculate the Riemann curvature tensor for three horizontal
vector fields at $(p, 0)$ we have
\begin{eqnarray}
\hat{R}^{f}(X^{h},Y^{h})Z^{h}&=&R(X,Y)Z+A_{f}(Y,Z)-A_{f}(X,Z)+A_{f}\Big(X,\nabla_{Y}Z+A_{f}(Y,Z)\Big)\nonumber\\
                 &&-A_{f}\Big(Y,\nabla_{X}Z+A_{f}(X,Z)\Big)-A_{f}([X,Y],Z)=0,
\end{eqnarray}
then $R=0$ and $f=C(constant)$.
\end{proof}
\begin{cor}\label{co:36}
Let $(M, g)$ be a Riemannian manifold and $TM$ be its tangent bundle with the rescaled Sasaki metric $\hat{g}^{f}$. If $f\neq C(constant)$, then
$(TM,\hat{g}^{f})$ is unflat.
\end{cor}

For the sectional curvatures of the tangent bundle we have the following.
\begin{prop}\label{pr: 37}
Let $(M, g)$ be a Riemannian manifold and equip the tangent bundle $(TM, \hat{g}^{f})$  with
the rescaled Sasaki metric $\hat{g}^{f}$. Let $(p,u)\in TM$ and $X,Y\in T_{p}M$ be two orthonormal tangent vectors at $p$. Let
$\hat{K}^{f}(X^{i},Y^{j})$ denote the sectional curvature of the plane spanned by $X^{i}$ and $Y^{j}$ with $i,j\in \{h,v\}$.
Then we have the following
\begin{eqnarray}
i) \ \hat{K}^{f}_{(p,u)}(X^{v},Y^{v})&=&0,\\
ii)\ \hat{K}^{f}_{(p,u)}(X^{h},Y^{v})&=&\frac{1}{4f^{2}(p)}|R(u,Y)X|^{2},\\
iii) \ \hat{K}^{f}_{(p,u)}(X^{h},Y^{h})&=&\frac{1}{f(p)}K(X,Y)-\frac{3}{4f^{2}(p)}|R(X,Y)u|^{2}+L_{f}(X,Y) ,
\end{eqnarray}
where
\begin{eqnarray*}
L_{f}(X,Y)&=&\frac{1}{f}\Big(g(\nabla_{X}A_{f}(Y,Y)-\nabla_{Y}A_{f}(X,Y),X)-g(A_{f}(X,\nabla_{Y}Y+A_{f}(Y,Y)),X)\nonumber\\
           &&-g(A_{f}(Y,\nabla_{X}Y+A_{f}(X,Y)),X)-g(A_{f}([X,Y],Y),X)\Big)     .
\end{eqnarray*}

\end{prop}
\begin{proof}
$i)$ It follows directly from Proposition 3.4 that the sectional curvature for a plane spanned by two vertical vectors vanishes.

$ii)$ Applying part $ii)$ of proposition 3.4 we get
\begin{eqnarray}
  \hat{K}^{f}(X^{h},Y^{v})&=&\frac{\hat{g}^{f}(\hat{R}^{f}(X^{h},Y^{v})Y^{v},X^{h})}
                                   {\hat{g}^{f}(X^{h},X^{h})\hat{g}^{f}(Y^{v},Y^{v})} \nonumber\\
                                  &=&\frac{1}{f}\Big(-\frac{1}{2f}\hat{g}^{f}((R(Y,Y)X)^{h},X^{h})
                                     -\frac{1}{4f^{2}}\frac{\hat{g}^{f}(R(u,Y)R(u,Y)X,X^{h})}{f g(R(u,Y)R(u,Y)X,X)}\Big) \nonumber\\
                                  &=&\frac{1}{4f^{2}}g(R(u,Y)X, R(u,Y)X) =\frac{1}{4f^{2}}|R(u,Y)X|^{2}.
\end{eqnarray}

$iii)$ It follows immediately from proposition 3.4 that
\begin{eqnarray}
  \hat{K}^{f}(X^{h},Y^{h})&=&\frac{1}{f^{2}}\hat{g}^{f}(\hat{R}^{f}(X^{h},Y^{h})Y^{h},X^{h}) \nonumber\\
                                  &=&\frac{1}{f}g(R(X,Y)Y,X)+\frac{3}{4f^{2}}g(R(Y,X)u,R(X,Y)u)\nonumber\\
                                  &=&\frac{1}{f}K(X,Y)-\frac{3}{4f^{2}}|R(X,Y)u|^{2}+\frac{1}{f}\Big(g(\nabla_{X}A_{f}(Y,Y)\nonumber\\
                                  &&-\nabla_{Y}A_{f}(X,Y),X)-g(A_{f}(X,\nabla_{Y}Y+A_{f}(Y,Y)),X)\nonumber\\
                                 &&-g(A_{f}(Y,\nabla_{X}Y+A_{f}(X,Y)),X)-g(A_{f}([X,Y],Y),X)\Big).
\end{eqnarray}
\end{proof}

\begin{thm}\label{th:38}
Let $(M, g)$ be a Riemannian manifold and equip the tangent bundle $(TM, \hat{g}^{f})$  with the rescaled Sasaki metric $\hat{g}^{f}$.
If the sectional curvature of $(TM,\hat{g}^{f})$ is upper bounded, then $(M,g)$ is flat; if $M$ compact and the sectional curvature
of $(TM,\hat{g}^{f})$ is lower bounded, then $(M,g)$ is flat.
\end{thm}
\begin{proof}
The statement follows directly from Proposition 3.7.
\end{proof}
\begin{prop}\label{pr: 39}
Let $(M, g)$ be a Riemannian manifold and equip the tangent bundle $(TM, \hat{g}^{f})$  with
the rescaled Sasaki metric $\hat{g}^{f}$. Let $(p,u)\in TM$ and $X,Y\in T_{p}M$ be two orthonormal tangent vectors at $p$. Let
$S$ denote the scalar curvature of $g$ and $\hat{S}^{f}$ denote the scalar curvature of $\hat{g}^{f}$. Then the following equation
holds
\begin{equation}
\hat{S}^{f}=\frac{1}{f}S-\frac{1}{4f^{2}}\sum_{i,j=1}^{m}\mid R(X_{i},Y_{j})u\mid^{2}+\sum_{i,j=1}^{m}L_{f}(X_{i},Y_{j})
\end{equation}
where $\{X_{1},\cdots,X_{m}\}$ is a local orthonormal frame for $TM$.
\end{prop}
\begin{proof}
For a local orthonormal frame $\{\frac{1}{\sqrt{f}}Y_{1},\cdots,\frac{1}{\sqrt{f}}Y_{m}, Y_{m+1},\cdots, Y_{2m}\}$ for $TTM$ with
$X^{h}_{i}=Y_{i}$ and $X^{v}_{i}=Y_{m+i}$ we get from proposition 3.7
\begin{eqnarray}
  \hat{S}^{f}&=&\sum_{i,j=1}^{m}\hat{K}^{f}(\frac{1}{\sqrt{f}}X^{h}_{i},\frac{1}{\sqrt{f}}X^{h}_{j})
           +2\sum_{i,j=1}^{m}\hat{K}^{f}(\frac{1}{\sqrt{f}}X^{h}_{i},X^{v}_{j})
           +\sum_{i,j=1}^{m}\hat{K}^{f}(X^{v}_{i},X^{v}_{j}) \nonumber\\
         &=&\sum_{i,j=1}^{m}[\hat{K}^{f}(X^{h}_{i},X^{h}_{j})+2\hat{K}^{f}(X^{h}_{i},X^{v}_{j})+\hat{K}^{f}(X^{v}_{i},X^{v}_{j})] \nonumber\\
         &=&\sum_{i,j=1}^{m}[\frac{1}{f}K(X_{i},X_{j})-\frac{3}{4f^{2}}\mid R(X_{i},X_{j})u\mid^{2}+L_{f}(X_{i},X_{j})]
             +2\sum_{i,j=1}^{m}\frac{1}{4f^{2}}|R(X_{j},u)X_{i}|^{2}.
\end{eqnarray}
In order to simplify this last expression we put $u=\sum_{i=1}^{m}u_{i}X_{i}$ we get

\begin{eqnarray}
\sum_{i,j=1}^{m}|R(X_{j},u)X_{i}|^{2}
                                &=&\sum_{i,j,k,l=1}^{m}u_{k}u_{l}g(R(X_{j},X_{k})X_{i},R(X_{j},X_{l})X_{i})\nonumber\\
                                &=&\sum_{i,j,k,l,s=1}^{m}u_{k}u_{l}g(R(X_{j},X_{k})X_{i},X_{s})g(R(X_{j},X_{l})X_{i},X_{s})\nonumber\\
                                &=&\sum_{i,j,k,l,s=1}^{m}u_{k}u_{l}g(R(X_{s},X_{i})X_{k},X_{j})g(R(X_{s},X_{i})X_{l},X_{j})\nonumber\\
                                &=&\sum_{i,j,k,l=1}^{m}u_{k}u_{l}g(R(X_{j},X_{i})X_{k},R(X_{j},X_{i})X_{l})\nonumber\\
                                &=&\sum_{i,j=1}^{m}|R(X_{j},X_{i})u|^{2}.
\end{eqnarray}
This completes the proof.
\end{proof}
\begin{cor}\label{co:40}
Let $(M, g)$ be a Riemannian manifold and $TM$ be its tangent bundle with the rescaled Sasaki metric $\hat{g}^{f}$.
Then $(TM,\hat{g}^{f})$ has constant scalar curvature if and only if $(M,g)$ is flat and $\sum_{i,j=1}^{m}L_{f}(X_{i},X_{j})=C(constant)$.
\end{cor}
\begin{proof}
The statement follows directly from Proposition 3.9.
\end{proof}

\section{Geodesics of The Rescaled Sasaki Metric}
Let $M$ be a Riemannian manifold with metric $g$. We denote by $\Im^{p}_{q}(M)$ the set of all tensor fields of type $(p,q)$ on
$M$. Manifolds, tensor fields and connections are always assumed to be differentiable and of class $C^{\infty}$. Let $T(M)$ be a
tangent bundle bundle of $M$, and $\pi$ the projection $\pi:T(M)\rightarrow M$. Let the manifold $M$ be covered
by system of coordinate neighbourhoods $(U,x^{i})$, where $(x^{i}), i=1,\cdots, n$ is a local coordinate system defined in the
neighbourhood $U$. Let $y^{i}$ be the Cartesian coordinates in each tangent spaces $T_{p}(M)$ and $P\in M$ with respect to the
natural base $\frac{\partial}{\partial x^{i}}$, $P$ being an arbitrary point in $U$ whose coordinates are $x^{i}$. Then we can
introduce local coordinates $(x^{i},y^{i})$ in open set $\pi^{-1}(U)\subset T(M_{n})$. We call them coordinates induced in
$\pi^{-1}(U)$ from $(U,x^{i})$. The projection $\pi$ is represented by $(x^{i},y^{i})\rightarrow (x^{i})$. The indices $i,j,\cdots$
 run from $1$ to $2n$.

Let $\hat{C}$ be a curve on $T(M_{n})$ and locally expressed by $x=x(\sigma),$ with respect to induced coordinates
$\frac{\partial}{\partial x_{i}}$ in $\pi^{-1}(U)\subset T(M_{n})$. The curve $\hat{C}$ is said to be a lift of
the curve $C$ and denoted by $C^{h}=(x(\sigma),x'(\sigma))$. The tangent vector field of $\hat{C}$ defined by
$T=(\frac{dx}{dt},\frac{dy}{dt})=x'^{h}+(\nabla_{x'}y)^{v}$. If the curve $\hat{C}$ is a geodesic, we get
$\nabla_{x'}x'=0$, then $y=x'$.

\begin{thm}\label{th:41}
Let $C$ be a geodesic on $T(M)$, if $f\neq c(constant)$ in any geodesics on $M$ , then the curve $C$
cannot be lifted to the geodesic of $\hat{g}^{f}$.
\end{thm}
\begin{proof}
By applying Proposition 3.2 we have
\begin{eqnarray}
 \hat{\nabla}^{f}_{x'^{h}+(\nabla_{x'}y)^{v}}(x'^{h}+(\nabla_{x'}y)^{v})
     &=&(\nabla_{x'}x')^{h}+A_{f}(x',x')^{h}-\frac{1}{2}(R(x',x')u)^{v} \nonumber\\
     &&+(\nabla_{x'}\nabla_{x'}y)^{v}+\frac{1}{2f}[R_{p}(u,\nabla_{x'}y)x']^{h}
     +\frac{1}{2f}[R_{p}(u,\nabla_{x'}y)x']^{h}  \nonumber\\
     &=&(\nabla_{x'}x')^{h}+\frac{1}{f}[R_{p}(u,\nabla_{x'}y)x']^{h}+A_{f}(x',x')^{h}+(\nabla_{x'}\nabla_{x'}y)^{v}.
\end{eqnarray}
For the curve $C$ is a geodesic on $M_{n}$, with respect to the adapted frame and taking account of
 $\hat{\nabla}^{f}_{T}T=0 $, then we get
\begin{eqnarray}
 && (a) \ \nabla_{x'}x'=-\frac{1}{f(x(t))}[R_{p}(y(t),\nabla_{x'}y(t))x'(t)]^{h}-A_{f}(x',x')^{h} ,  \nonumber\\
 && (b) \  \nabla_{x'}\nabla_{x'}y=0.
\end{eqnarray}
Applying part $iv)$ of proposition 3.4 and
\begin{equation}
A_{f}(x',x')=\frac{1}{2f}[2x'(t)x'-g(x',x')\texttt{grad} f],
\end{equation}
if $\langle A_{f}(x',x'), x'\rangle =0$, we get
\begin{equation}
2X'(f)g(x',x')-g(x',x')X'(f)=g(x',x')X'(f)=0.
\end{equation}
Then $X'(f)=0$, $grad(f)=0$ and
$\frac{\texttt{d}f(x(t))}{\texttt{d}t}=0$, so we get $f=c(constant)$ in any
geodesics on $M$.
\end{proof}

\begin{cor}\label{co:42}
If $(x(t),y(t))$ is geodesic and $|y(t)|=C$, then $\nabla_{x'}x'=-A_{f}(x',x')$.
\end{cor}
\begin{proof}
By applying $(a)$ of equation $(4.2)$ we have
\begin{equation}
0=\nabla_{x'}\langle y, y\rangle =\langle \nabla_{x'}y, y\rangle
+\langle y, \nabla_{x'}y\rangle,
\end{equation}
and
\begin{equation}
0=\nabla_{x'}\langle \nabla_{x'}y, y\rangle =\langle
\nabla_{x'}\nabla_{x'}y, y\rangle +\langle \nabla_{x'}y,
\nabla_{x'}y\rangle.
\end{equation}
Then we get $\langle \nabla_{x'}y, y\rangle=0$ and $\nabla_{x'}y=0$,
from which the result follows.
\end{proof}

\begin{thm}\label{th:43}
Let $C_{1}$ and $C_{2}$ be two geodesics on $M_{n}$ departure from the same arbitrary point,
and their initial tangent vectors are not parallel. If the lifts of two geodesics on $M$ are
geodesics on $T(M)$ with the metric $\hat{g}^{f}$, then
$f=c(constant)$.
\end{thm}
\begin{proof}
By applying $(a)$ of equation $(4.2)$ we have
\begin{equation}
2X'(f)x'-g(x',x')X'(f)=2\tilde{X'}(f)\tilde{x}'-g(\tilde{x}',\tilde{x}')X'(f).
\end{equation}
Using $X'(0) \nparallel\tilde{X'}(0)$ we get
$\texttt{grad}f(x_{0})=0$, then we obtain $f=c(constant)$.
\end{proof}

The submersion geodesic $C$ is said to be the image under $\pi$ of the geodesic $\hat{C}$ on $TM$.
Let $C=\pi\circ \hat{C}$ be a submersion geodesic on $M$, then $\hat{\nabla}^{f}_{T}T=0 $.
Using this condition we have
\begin{thm}\label{th:44}
Let $M$ be a flat manifold, the submersion geodesic is always geodesics on $M$, then $f=c(constant)$.
\end{thm}

\section{The Rescaled Cheeger-Gromoll Metric}
In \cite{CG}, Cheeger and Gromoll studied complete manifolds of nonnegative curvature and suggest a construction of Riemannian metrics useful
 in that context. This can be used to obtain a natural metric $\tilde{g}^{f}$ on the tangent bundle $TM$ of a given Riemannian manifold $(M, g)$.

 For a vector field $u\in C^{\infty}(TM)$ we shall by $U$ denote its canonical vertical vector field on $TM$ which in local coordinates
 is given by
\begin{equation}
U=\sum_{i=1}^{m}v_{m+i}(\frac{\partial}{\partial v_{m+i}})_{(p,u)},
\end{equation}
where $u=(v_{m+1},\cdots,v_{2})$. To simplify our notation we define the function $r: TM\rightarrow\mathbb{R}$
by $r(p,u)=|u|=\sqrt{g_{p}(u,u)}$
and $\alpha=1+r^{2}.$
\begin{defn}
Let $(M,g)$ be a Riemannian manifold. Let $f>0$ and $f\in C^{\infty}(M)$. Then the rescaled Cheeger-Gromoll metric $\tilde{g}^{f}$ on the tangent
bundle $TM$ of $M$ is given by
\begin{eqnarray}
i) \ \tilde{g}_{(p,u)}^{f}(X^{h}, Y^{h})&=&f(p)g_{p}(X, Y),\\
ii) \ \tilde{g}_{(p,u)}^{f}(X^{v},  Y^{h})&=&0 ,\\
iii) \ \tilde{g}_{(p,u)}^{f}(X^{v},  Y^{v})&=&\frac{1}{1+r^{2}}(g_{p}(X, Y)+g_{p}(X, u)g_{p}(Y, u))
\end{eqnarray}
for all vector fields $X, Y\in C^{\infty}(TM)$.
\end{defn}
It is obvious that the rescaled Cheeger-Gromoll metric $\tilde{g}^{f}$ is contained in the class of rescaled natural metrics introduced earlier.

\begin{prop}\label{pr: 42}
Let $(M, g)$ be a Riemannian manifold and $\tilde{\nabla}^{f}$ be Levi-Civita connection of the tangent bundle $(TM, \tilde{g}^{f})$ equipped
with the rescaled Cheeger-Gromoll metric. Then
\begin{eqnarray}
i) \ (\tilde{\nabla}^{f}_{X^{h}}Y^{h})_{(p,u)}&=&(\nabla_{X}Y)^{h}_{(p,u)}+\frac{1}{2f(p)}
   \Big((X(f)Y+Y(f)X)-g(X,Y)(\texttt{d}f)^{*}\Big)^{h}_{p}\nonumber\\
&&-\frac{1}{2}\Big(R_{p}(X,Y)u\Big)^{v},\\
ii) \ (\tilde{\nabla}^{f}_{X^{h}}Y^{v})_{(p,u)}&=&(\nabla_{X}Y)^{v}_{(p,u)}+\frac{1}{2\alpha f(p)}\Big(R_{p}(u,Y)X\Big)^{h},\\
iii) \ (\tilde{\nabla}^{f}_{X^{v}}Y^{h})_{(p,u)}&=&\frac{1}{2\alpha f(p)}\Big(R_{p}(u,X)Y\Big)^{h},\\
iv) \ (\tilde{\nabla}^{f}_{X^{v}}Y^{v})_{(p,u)}&=&-\frac{1}{\alpha}\Big(\tilde{g}_{(p,u)}^{f}(X^{v}, U)Y^{v}
   +\tilde{g}_{(p,u)}^{f}(Y^{v}, U)X^{v}\Big)\nonumber\\
&&+\frac{1+\alpha}{\alpha}\tilde{g}_{(p,u)}^{f}(X^{v},Y^{v})U-\frac{1}{\alpha}\tilde{g}_{(p,u)}^{f}(X^{v}, U)\tilde{g}_{(p,u)}^{f}(Y^{v}, U)U
\end{eqnarray}
for any $X,Y\in C^{\infty}(TM)$, $\xi=(p, u)\in TM$.
\end{prop}
\begin{proof}
 $i)$ The statement is a direct consequence of Corollary 2.3.

 $ii)$ By applying Lemma 2.2 and Definition 4.1 we get
\begin{eqnarray}
  2\tilde{g}(\tilde{\nabla}^{f}_{X^{h}}Y^{v}, Z^{h})&=&-\frac{1}{f}\tilde{g}(Y^{v},(R(Z,X)u)^{v})\nonumber\\
                                              &=&-\frac{1}{\alpha f}\Big(g(Y,R(Z,X)u)+g(Y,u)g(R(Z,X)u,u)\Big)\nonumber\\
                                              &=&\frac{1}{\alpha f}g\Big(R(u,Y)X,Z\Big).
\end{eqnarray}
From Definition 3.7 and Lemma 4.1 in \cite{GK2} it follows that
\begin{equation}
X^{h}(\frac{1}{\alpha})=0 \quad and \quad X^{h}(g(Y,u))\circ\pi=g(\nabla_{X}Y,u)\circ\pi,
\end{equation}
so
\begin{equation}
X^{h}(\tilde{g}^{f}(Y^{v},Z^{v}))=\tilde{g}^{f}((\nabla_{X}Y)^{v},Z^{v})+\tilde{g}^{f}(Y^{v},(\nabla_{X}Z)^{v}).
\end{equation}
This means that
\begin{eqnarray}
  2\tilde{g}^{f}(\tilde{\nabla}^{f}_{X^{h}}Y^{v}, Z^{v})&=&X^{h}(\tilde{g}^{f}(Y^{v},Z^{v}))+\tilde{g}^{f}(Z^{v},(\nabla_{X}Y)^{v})
                                                          -\tilde{g}^{f}(Y^{v},(\nabla_{X}Z)^{v})  \nonumber\\
                                                      &=&2\tilde{g}^{f}((\nabla_{X}Y)^{v},Z^{v}).
\end{eqnarray}

$iii)$ Calculations similar to those in $ii)$ give
\begin{eqnarray}
  2\tilde{g}(\tilde{\nabla}^{f}_{X^{v}}Y^{h}, Z^{h})&=&\frac{1}{f}\tilde{g}(X^{v},(R(Y,Z)u)^{v})
                                                         =\frac{1}{\alpha f}\tilde{g}((R(u,X)Y)^{h},Z^{h}).
\end{eqnarray}
The rest follows by
\begin{eqnarray}
  2\tilde{g}^{f}(\tilde{\nabla}^{f}_{X^{v}}Y^{h}, Z^{v})&=&Y^{h}(\tilde{g}^{f}(Z^{v},X^{v}))-\tilde{g}^{f}(Z^{v},(\nabla_{Y}X)^{v})
                                                          -\tilde{g}^{f}(X^{v},(\nabla_{Y}Z)^{v})  \nonumber\\
                                              &=&\tilde{g}^{f}(Z^{v},(\nabla_{Y}X)^{v})+\tilde{g}^{f}(X^{v},(\nabla_{Y}Z)^{v})
                                              -\tilde{g}^{f}(Z^{v},(\nabla_{Y}X)^{v})-\tilde{g}^{f}(X^{v},(\nabla_{Y}Z)^{v})\nonumber\\
                                              &=&0.
\end{eqnarray}

$iv)$ Using Lemma 2.2 we yield
\begin{eqnarray}
  2f\tilde{g}^{f}(\tilde{\nabla}^{f}_{X^{v}}Y^{v}, Z^{h})&=&-Z^{h}(\tilde{g}^{f}(X^{v},Y^{v}))+\tilde{g}^{f}(Y^{v},(\nabla_{Z}X)^{v})
                                                         +\tilde{g}^{f}(X^{v},(\nabla_{Z}Y)^{v})  \nonumber\\
                                              &=&-\tilde{g}^{f}(Y^{v},(\nabla_{Z}X)^{v})-\tilde{g}^{f}(X^{v},(\nabla_{Z}Y)^{v})
                                              +\tilde{g}^{f}(Y^{v},(\nabla_{Z}X)^{v})+\tilde{g}^{f}(X^{v},(\nabla_{Z}Y)^{v}) \nonumber\\
                                              &=&0.
\end{eqnarray}
Using $X^{v}(f(r^{2}))=2f'(r^{2})g(X,u)$ and $\alpha=1+r^{2}$ we get
\begin{eqnarray}
  X^{v}\tilde{g}^{f}(Y^{v}, Z^{v})&=&-\frac{2}{\alpha^{2}}g(X,u)\Big(g(Y,Z)+g(Y,u)g(Z,u)\Big)  \nonumber\\
                                              &&+\frac{1}{\alpha}\Big(g(X,Y)g(Z,u)+g(X,Z)g(Y,u)\Big).
\end{eqnarray}
The definition of the rescaled Cheeger-Gromoll metric implies that
\begin{equation}
\tilde{g}^{f}(X^{v}, U)=\frac{1}{\alpha}\Big(g(X,u)+g(X,u)g(u,u)\Big)=g(X,u).
\end{equation}
This leads to the following
\begin{eqnarray}
 \alpha^{2}\tilde{g}^{f}(\tilde{\nabla}^{f}_{X^{v}}Y^{v}, Z^{v})&=&\frac{\alpha^{2}}{2}\Big(X^{v}(\tilde{g}^{f}(Y^{v},Z^{v}))
                        +Y^{v}(\tilde{g}^{f}(Z^{v},X^{v}))-Z^{v}(\tilde{g}^{f}(X^{v},Y^{v}))\Big)  \nonumber\\
                                         &=&-g(X,u)\Big(g(Y,Z)+g(Y,u)g(Z,u)\Big) \nonumber\\
                                         &&+\frac{\alpha}{2}\Big(g(X,Y)g(Z,u)+g(X,Z)g(Y,u)\Big)\nonumber\\
                                         &&-g(Y,u)\Big(g(Z,X)+g(Z,u)g(X,u)\Big) \nonumber\\
                                         &&+\frac{\alpha}{2}\Big(g(Y,Z)g(X,u)+g(Y,X)g(Z,u)\Big) \nonumber\\
                                         &&+g(Z,u)\Big(g(X,Y)+g(X,u)g(Y,u)\Big) \nonumber\\
                                         &&-\frac{\alpha}{2}\Big(g(Z,X)g(Y,u)+g(Z,Y)g(X,u)\Big) \nonumber\\
                                         &=&g\Big(\big(g(X,Y)-g(X,u)g(Y,u)\big)u+\alpha g(X,Y)u \nonumber\\
                                         &&-g(X,u)Y-g(Y,u)X,Z\Big).
\end{eqnarray}
By using the definition of the metric we see that this gives the statement to proof.
\end{proof}
Having determined the Levi-Civita connection we are ready to calculate the Riemann curvature tensor of $TM$. But first we state the
following useful Lemma.
\begin{lem}\label{le:43}
Let $(M, g)$ be a Riemannian manifold and $\tilde{\nabla}^{f}$ be the Levi-Civita connection of the tangent bundle $(TM, \tilde{g}^{f})$,
equipped with the rescaled Cheeger-Gromoll metric $\tilde{g}^{f}$. Let $F:TM\rightarrow TM$ is a smooth bundle endomorphism of the
tangent bundle, then
\begin{eqnarray}
  (\tilde{\nabla}^{f}_{X^{v}}F^{v})_{\xi}&=&F(X)^{v}_{\xi}-\frac{1}{\alpha}\Big(\tilde{g}^{f}(X^{v},U)F^{v}+\tilde{g}^{f}(F^{v},U)X^{v} \nonumber\\
                                              &&-(1+\alpha)\tilde{g}^{f}(F^{v},X^{v})U+\tilde{g}^{f}(X^{v},U)\tilde{g}^{f}(F^{v},U)U\Big)_{\xi}
\end{eqnarray}
and
\begin{equation}
(\tilde{\nabla}^{f}_{X^{v}}F^{h})_{\xi}=F(X)^{h}_{\xi}+\frac{1}{2\alpha f(p)}\Big(R(u,X)F(u)\Big)^{h}_{\xi}
\end{equation}
for any $X\in C^{\infty}(TM)$ and $\xi=(p, u)\in TM$.
\end{lem}
\begin{proof}
The statement is a direct consequence of Lemma 2.5 and Proposition 5.2.
\end{proof}
\begin{prop}\label{pr: 54}
Let $(M, g)$ be a Riemannian manifold and $\tilde{R}^{f}$ be the Riemann curvature tensor of the tangent bundle $(TM, \tilde{g}^{f})$
equipped with the rescaled Sasaki metric. Then the following formulae hold
\begin{eqnarray*}
i) \ \tilde{R}^{f}(X^{h},Y^{h})Z^{h}&=&\nabla_{X}(\nabla_{Y}Z+A_{f}(Y,Z))^{h}+A_{f}(X,\nabla_{Y}Z+A_{f}(Y,Z))^{h} \nonumber\\
                            &&-\frac{1}{2}[R(X,\nabla_{Y}Z+A_{f}(Y,Z))u]^{v}-\nabla_{Y}(\nabla_{X}Z+A_{f}(X,Z))^{h}\nonumber\\
                            &&-A_{f}(Y,\nabla_{X}Z+A_{f}(X,Z))^{h}+\frac{1}{2}\Big(R(Y,\nabla_{X}Z+A_{f}(X,Z))u\Big)^{v}\nonumber\\
                            &&-(\nabla_{[X,Y]}Z)^{h}-A_{f}([X,Y],Z)^{h}-\frac{1}{2}\Big(R([X,Y],Z)u\Big)^{v}\nonumber\\
                             &&+\frac{1}{2\alpha f}\Big(R(u,R(X,Y)u)Z\Big)^{h}\nonumber\\
\end{eqnarray*}
\begin{eqnarray}
                               &&+\frac{1}{2}[\nabla_{Y}(R(X,Z)u)]^{v}+\frac{1}{4\alpha f(p)}\Big(R(u,R(X,Z)u)Y\Big)^{h} \nonumber\\
                               &&-\frac{1}{2}[\nabla_{X}(R(Y,Z)u)]^{v}-\frac{1}{4\alpha f(p)}\Big(R(u,R(Y,Z)u)X\Big)^{h},
\end{eqnarray}
\begin{eqnarray}
ii) \ \tilde{R}^{f}(X^{h},Y^{h})Z^{v}&=&(R(X,Y)Z)^{v}+\frac{1}{2\alpha }\Big(\nabla_{Z}(\frac{1}{f}R(u,Z)Y)
                                 -\nabla_{Y}(\frac{1}{f}R(u,Z)X)\Big)^{h} \nonumber\\
                                           &&-\frac{1}{4\alpha f(p)}\Big(R(X,R(u,Z)Y)u-R(Y,R(u,Z)X)u\Big)^{v} \nonumber\\
                                           &&+\frac{1}{\alpha}\Big(A_{f}(X,\frac{1}{2f}R(u,Z)Y)-A_{f}(Y,\frac{1}{2f}R(u,Z)X)\Big)^{h}  \nonumber\\
                  &&-\frac{1}{\alpha}\tilde{g}^{f}(Z^{v},u)(R(X,Y)u)^{v}+\frac{1+\alpha}{\alpha}\tilde{g}^{f}((R(X,Y)u)^{v},Z^{v})U.
\end{eqnarray}
\begin{eqnarray}
iii) \ \tilde{R}^{f}(X^{h},Y^{v})Z^{h} &=& \frac{1}{2\alpha }\tilde{\nabla}^{f}_{X^{h}}(\frac{1}{f}R(u,Y)Z)^{h}\nonumber\\
                         &&-\frac{1}{2\alpha f}(R(u,\nabla_{X}Y)Z)^{h}-\frac{1}{2\alpha f}(R(u,Y)\nabla_{X}Z)^{h}
                           +\frac{1}{2}(R(X,Z)Y)^{v}  \nonumber\\
                     &&-\frac{1}{2\alpha}\tilde{g}^{f}(Y^{v},U)(R(X,Z)u)^{v}-\frac{1}{2\alpha}\tilde{g}^{f}((R(X,Z)u)^{v},U)Y^{v}\nonumber\\
                     &&+\frac{1+\alpha}{2\alpha}\tilde{g}^{f}((R(X,Z)u)^{v},Y^{v})U
                        -\frac{1}{2\alpha}\tilde{g}^{f}(Y^{v},U)\tilde{g}^{f}((R(X,Z)u)^{v},U)U \nonumber\\
                     && -\frac{1}{2\alpha f}(R(u,Y)A_{f}(X,Z))^{h},
\end{eqnarray}
\begin{eqnarray}
iv)\ \tilde{R}^{f}(X^{h},Y^{v})Z^{v}&=&-\frac{1}{2\alpha f}\Big(R(Y,Z)X\Big)^{h}
                                -\frac{1}{4\alpha^{2} f^{2}}\Big(R(u,Y)R(u,Z)X\Big)^{h}\nonumber\\
                               &&+\frac{1}{2\alpha^{2} f}[g(Y,U)(R(u,Z)X)^{h}-g(Z,u)(R(u,Y)X)^{h}],
\end{eqnarray}
\begin{eqnarray}
v)\ \tilde{R}^{f}_{(p,u)}(X^{v},Y^{v})Z^{h}&=&-\frac{1}{2\alpha f}\Big(R(X,Y)Z\Big)^{h}
                                            -\frac{1}{4\alpha^{2} f^{2}}\Big(R(u,X)R(u,Y)Z\Big)^{h}\nonumber\\
                                  &&+\frac{1}{2\alpha f}\Big(R(Y,X)Z\Big)^{h}
                                            +\frac{1}{4\alpha^{2} f^{2}}\Big(R(u,Y)R(u,X)Z\Big)^{h},
\end{eqnarray}
\begin{eqnarray}
vi) \ \tilde{R}^{f}_{(p,u)}(X^{v},Y^{v})Z^{v}&=&\frac{1+\alpha+\alpha^{2}}{\alpha^{2}}
                                      (\tilde{g}^{f}(Y^{v},Z^{v})X^{v}-\tilde{g}^{f}(X^{v},Z^{v})Y^{v})+ \nonumber\\
                               &&+\frac{2+\alpha}{\alpha^{2}}
                                      (\tilde{g}^{f}(X^{v},Z^{v})g(Y,u)U-\tilde{g}^{f}(Y^{v},Z^{v})g(X,u)U)+ \nonumber\\
                               &&+\frac{2+\alpha}{\alpha^{2}}
                                      (g(X,u)g(Z,u)Y^{v}-g(Y,u)g(Z,u)X^{v}).
\end{eqnarray}
for any $X, Y, Z\in T_{p}M$.
\end{prop}
\begin{proof}
$i)$ By $i)$ of Proposition 4.2 and direct calculation we get
\begin{eqnarray}
\tilde{R}^{f}(X^{h},Y^{h})Z^{h}&=&\tilde{\nabla}^{f}_{X^{h}}\tilde{\nabla}^{f}_{Y^{h}}Z^{h}-\tilde{\nabla}^{f}_{Y^{h}}\tilde{\nabla}^{f}_{X^{h}}Z^{h}
                                 -\tilde{\nabla}^{f}_{[X^{h},Y^{h}]}Z^{h}\nonumber\\
                             &=&\tilde{\nabla}^{f}_{X^{h}}((\nabla_{Y}Z)^{h}+A_{f}(Y,Z)^{h}-\frac{1}{2}(R(Y,Z)u)^{v})  \nonumber\\
                            && -\tilde{\nabla}^{f}_{Y^{h}}((\nabla_{X}Z)^{h}+A_{f}(X,Z)^{h}-\frac{1}{2}(R(X,Z)u)^{v})\nonumber\\
                            &&-\tilde{\nabla}^{f}_{[X,Y]^{h}-(R(X,Y)u)^{v}}Z^{h}\nonumber\\
                            &=&\nabla_{X}(\nabla_{Y}Z+A_{f}(Y,Z))^{h}+A_{f}(X,\nabla_{Y}Z+A_{f}(Y,Z))^{h} \nonumber\\
                            &&-\frac{1}{2}\Big(R(X,\nabla_{Y}Z+A_{f}(Y,Z))u\Big)^{v}-\nabla_{Y}(\nabla_{X}Z+A_{f}(X,Z))^{h}\nonumber\\
                            &&-A_{f}(Y,\nabla_{X}Z+A_{f}(X,Z))^{h}+\frac{1}{2}\Big(R(Y,\nabla_{X}Z+A_{f}(X,Z))u\Big)^{v}\nonumber\\
                            &&-(\nabla_{[X,Y]}Z)^{h}-A_{f}([X,Y],Z)^{h}-\frac{1}{2}\Big(R([X,Y],Z)u\Big)^{v}\nonumber\\
                             &&+\frac{1}{2\alpha f}\Big(R(u,R(X,Y)u)Z\Big)^{h}\nonumber\\
                               &&+\frac{1}{2}\Big(\nabla_{Y}(R(X,Z)u)\Big)^{v}+\frac{1}{4\alpha f}\Big(R(u,R(X,Z)u)Y\Big)^{h} \nonumber\\
                               &&-\frac{1}{2}\Big(\nabla_{X}(R(Y,Z)u)\Big)^{v}-\frac{1}{4\alpha f}\Big(R(u,R(Y,Z)u)X\Big)^{h}.
\end{eqnarray}

$ii)$ Note that the equation $\tilde{g}^{f}_{(p,u)}(X^{v},U)=g_{p}(X,u)$ implies that
\begin{equation}
\tilde{g}^{f}_{(p,u)}((R(X,Y)u)^{v},U)=g_{p}(R(X,Y)u,u)=0,
\end{equation}
Hence
\begin{eqnarray*}
\alpha\tilde{R}^{f}(X^{h},Y^{h})Z^{v}&=&\alpha\tilde{\nabla}^{f}_{X^{h}}\tilde{\nabla}^{f}_{Y^{h}}Z^{v}
                                        -\alpha\tilde{\nabla}^{f}_{Y^{h}}\tilde{\nabla}^{f}_{X^{h}}Z^{v}
                                        -\alpha\tilde{\nabla}^{f}_{[X^{h},Y^{h}]}Z^{v}\nonumber\\
                        &=&\tilde{\nabla}^{f}_{X^{h}}(\alpha(\nabla_{Y}Z)^{v}+\frac{1}{2f}(R(u,Z)Y)^{h})  \nonumber\\
                            && -\tilde{\nabla}^{f}_{Y^{h}}(\alpha(\nabla_{X}Z)^{h}+\frac{1}{2f}(R(u,Z)X)^{h})
                               -\alpha\tilde{\nabla}^{f}_{[X,Y]^{h}-(R(X,Y)u)^{v}}Z^{v} \nonumber\\
&=&(\nabla_{X}(\frac{1}{2f}R(u,Z)Y))^{h}-\frac{1}{4f}(R(X,R(u,Z)Y)u)^{v}\nonumber\\
                             &&+\frac{1}{2f}(R(u,\nabla_{Y}Z)X)^{h}\nonumber\\
                               && +A_{f}(X,\frac{1}{2f}R(u,Z)Y)^{h}+\alpha(\nabla_{X}\nabla_{Y}Z)^{v} \nonumber\\
                            &&-(\nabla_{Y}(\frac{1}{2f}R(u,Z)X))^{h}+\frac{1}{4f}(R(Y,R(u,Z)X)u)^{v}\nonumber\\
                        &&-\frac{1}{2f}(R(u,\nabla_{X}Z)Y)^{h}\nonumber\\
                            &&-A_{f}(Y,\frac{1}{2f}R(u,Z)X)^{h}-\alpha(\nabla_{Y}\nabla_{X}Z)^{v} \nonumber\\
                            &&-\frac{1}{2f}(R(u,Z)[X,Y])^{h}-\alpha(\nabla_{[X,Y]}Z)^{v} \nonumber\\
\end{eqnarray*}
\begin{eqnarray}
                            &&-[\tilde{g}^{f}((R(X,Y)u)^{v},U)Z^{v}+\tilde{g}^{f}(Z^{v},U)(R(X,Y)u)^{v}] \nonumber\\
                            &&+(1+\alpha)\tilde{g}^{f}((R(X,Y)u)^{v},Z^{v})U-\tilde{g}^{f}((R(X,Y)u)^{v},U)\tilde{g}^{f}(Z^{v},U)U  \nonumber\\
                        &=&\alpha(R(X,Y)Z)^{v}+\frac{1}{2f}[\nabla_{Z}(R(u,Z)Y)-\nabla_{Y}(R(u,Z)X)]^{h} \nonumber\\
                                           &&-\frac{1}{4f}[R(X,R(u,Z)Y)u-R(Y,R(u,Z)X)u]^{v} \nonumber\\
                                           &&+[A_{f}(X,\frac{1}{2f}R(u,Z)Y)-A_{f}(Y,\frac{1}{2f}R(u,Z)X)]^{h}  \nonumber\\
                        &&-\tilde{g}^{f}(Z^{v},u)(R(X,Y)u)^{v}+(1+\alpha)\tilde{g}^{f}((R(X,Y)u)^{v},Z^{v})U.
\end{eqnarray}
$iii)$ Calculations similar to those above produce the third formula
\begin{eqnarray}
\tilde{R}^{f}(X^{h},Y^{v})Z^{h}&=&\tilde{\nabla}^{f}_{X^{h}}\tilde{\nabla}^{f}_{Y^{v}}Z^{h}
                                        -\tilde{\nabla}^{f}_{Y^{v}}\tilde{\nabla}^{f}_{X^{h}}Z^{h}
                                        -\tilde{\nabla}^{f}_{[X^{h},Y^{v}]}Z^{h}\nonumber\\
                        &=&\frac{1}{2\alpha }\tilde{\nabla}^{f}_{X^{h}}(\frac{1}{f}R(u,Y)Z)^{h}
                         -\tilde{\nabla}^{f}_{(\nabla_{X}Y)^{v}}Z^{h}  \nonumber\\
                            && -\tilde{\nabla}^{f}_{Y^{v}}[(\nabla_{X}Z)^{h}-\frac{1}{2}(R(X,Z)u)^{v}+A_{f}(X,Z)^{h}] \nonumber\\
                       &=& \frac{1}{2\alpha }\tilde{\nabla}^{f}_{X^{h}}(\frac{1}{f}R(u,Y)Z)^{h}\nonumber\\
                         &&-\frac{1}{2\alpha f}(R(u,\nabla_{X}Y)Z)^{h}-\frac{1}{2\alpha f}(R(u,Y)\nabla_{X}Z)^{h}
                           +\frac{1}{2}(R(X,Z)Y)^{v}  \nonumber\\
                     &&-\frac{1}{2\alpha}\tilde{g}^{f}(Y^{v},U)(R(X,Z)u)^{v}-\frac{1}{2\alpha}\tilde{g}^{f}((R(X,Z)u)^{v},U)Y^{v}\nonumber\\
                     &&+\frac{1+\alpha}{2\alpha}\tilde{g}^{f}((R(X,Z)u)^{v},Y^{v})U
                        -\frac{1}{2\alpha}\tilde{g}^{f}(Y^{v},U)\tilde{g}^{f}((R(X,Z)u)^{v},U)U \nonumber\\
                     && -\frac{1}{2\alpha f}(R(u,Y)A_{f}(X,Z))^{h}.
\end{eqnarray}
$iv)$ Since $X^{v}_{(p,u)}(f(r^{2}))=2f'(r^{2})g_{p}(X,u)$ and $(\tilde{\nabla}^{f}_{X^{h}}U)_{(p,u)}=0$ we get
\begin{eqnarray*}
2\alpha\tilde{R}^{f}(X^{h},Y^{v})Z^{v}&=&2\alpha[\tilde{\nabla}^{f}_{X^{h}}\tilde{\nabla}^{f}_{Y^{v}}Z^{v}
                                        -\tilde{\nabla}^{f}_{Y^{v}}\tilde{\nabla}^{f}_{X^{h}}Z^{v}
                                        -\tilde{\nabla}^{f}_{[X^{h},Y^{v}]}Z^{v}]\nonumber\\
                        &=&-2\tilde{\nabla}^{f}_{X^{h}}[\tilde{g}^{f}(Y^{v},U)Z^{v}
                        -(1+\alpha)\tilde{g}^{f}(Y^{v},Z^{v})U \nonumber\\
                        && +\tilde{g}^{f}(Z^{v},U)Y^{v} +\tilde{g}^{f}(Y^{v},U)\tilde{g}^{f}(Z^{v},U)U  ]  \nonumber\\
                        && -\alpha\tilde{\nabla}^{f}_{Y^{v}}(\frac{1}{\alpha f}R(u,Z)X)^{h}
                         -2\alpha[\tilde{\nabla}^{f}_{Y^{v}}(\nabla_{X}Z)^{v}+\tilde{\nabla}^{f}_{(\nabla_{X}Y)^{v}}Z^{v}] \nonumber\\
                      &=&-g(Y,u)[\frac{1}{\alpha f}(R(u,Z)X)^{h}+2(\nabla_{X}Z)^{v}]\nonumber\\
                      &&-g(Z,u)[\frac{1}{\alpha f}(R(u,Y)X)^{h}+2(\nabla_{X}Y)^{v}]\nonumber\\
              &&+\frac{2}{\alpha f}g(Y,u)(R(u,Z)X)^{h}\nonumber\\
                            &&-\tilde{\nabla}^{f}_{Y^{v}}(\frac{1}{f}R(u,Z)X)^{h}\nonumber\\
                      &&+2[g(Y,u)(\nabla_{X}Z)^{v}+g(\nabla_{X}Z,u)Y^{v}]\nonumber\\
                      &&-(1+\alpha)\tilde{g}^{f}(Y^{v},(\nabla_{X}Z)^{v})U+g(Y,u)g(\nabla_{X}Z,u)U\nonumber\\
\end{eqnarray*}
\begin{eqnarray}
                      &&+g(\nabla_{X}Y,u)Z^{v}+g(Z,u)(\nabla_{X}Y)^{v}\nonumber\\
                      &&-(1+\alpha)\tilde{g}^{f}((\nabla_{X}Y)^{v},Z^{v})U+g(\nabla_{X}Y,u)g(Z,u)U\nonumber\\
                      &=&-\tilde{\nabla}^{f}_{Y^{v}}(\frac{1}{f}R(u,Z)X)^{h}\nonumber\\
                               &&+\frac{1}{\alpha f}[g(Y,U)(R(u,Z)X)^{h}-g(Z,u)(R(u,Y)X)^{h}]\nonumber\\
                                &=&-\frac{1}{f}\Big(R(Y,Z)X\Big)^{h}
                                -\frac{1}{2\alpha f^{2}}\Big(R(u,Y)R(u,Z)X\Big)^{h}\nonumber\\
                               &&+\frac{1}{\alpha f}[g(Y,U)(R(u,Z)X)^{h}-g(Z,u)(R(u,Y)X)^{h}]
\end{eqnarray}
For the last equation we have to show that all the terms not containing the Riemann curvature tenson $R$ vanish. But since
\begin{equation}
\tilde{g}^{f}(Y^{v},(\nabla_{X}Z)^{v})U=\frac{1}{\alpha}[g(Y,\nabla_{X}Z)+g(Y,u)g(\nabla_{X}Z,u)]U,
\end{equation}
the rest becomes
\begin{equation}
-\frac{2}{\alpha}[g(Y,\nabla_{X}Z)+g(Y,u)g(\nabla_{X}Z,u)+g(Z,\nabla_{X}Y)+g(Z,u)g(\nabla_{X}Y,u)]U,
\end{equation}
which vanishes, because
\begin{equation}
-\frac{2}{\alpha}X^{h}[\tilde{g}^{f}(Y^{v},Z^{v})+\tilde{g}^{f}(Y^{v},U)\tilde{g}^{f}(Z^{v},U)]U=0.
\end{equation}
$v)$ First we notice that
\begin{eqnarray}
\tilde{\nabla}^{f}_{X^{v}}\tilde{\nabla}^{f}_{Y^{v}}Z^{h}&=&\frac{1}{2\alpha}\tilde{\nabla}^{f}_{X^{v}}(\frac{1}{f}R(u,Y)Z)^{h}\nonumber\\
                                                    &=&-\frac{1}{2\alpha f}\Big(R(X,Y)Z\Big)^{h}
                                            -\frac{1}{4\alpha^{2} f^{2}}\Big(R(u,X)R(u,Y)Z\Big)^{h}.
\end{eqnarray}

By using the fact that $[X^{v},Y^{v}]=0$ we get
\begin{eqnarray}
\tilde{R}^{f}(X^{v},Y^{v})Z^{h}&=&\tilde{\nabla}^{f}_{X^{v}}\tilde{\nabla}^{f}_{Y^{v}}Z^{h}
                                        -\tilde{\nabla}^{f}_{Y^{v}}\tilde{\nabla}^{f}_{X^{v}}Z^{h} \nonumber\\
                               &=&\frac{1}{2\alpha}\tilde{\nabla}^{f}_{X^{v}}(\frac{1}{f}R(u,Y)Z)^{h}
                                  -\frac{1}{2\alpha}\tilde{\nabla}^{f}_{Y^{v}}(\frac{1}{f}R(u,X)Z)^{h}\nonumber\\
                                  &=&-\frac{1}{2\alpha f}\Big(R(X,Y)Z\Big)^{h}
                                            -\frac{1}{4\alpha^{2} f^{2}}\Big(R(u,X)R(u,Y)Z\Big)^{h}\nonumber\\
                                  &&+\frac{1}{2\alpha f}\Big(R(Y,X)Z\Big)^{h}
                                            +\frac{1}{4\alpha^{2} f^{2}}\Big(R(u,Y)R(u,X)Z\Big)^{h}.
\end{eqnarray}
$vi)$ The result similar to Proposition 8.5 in \cite{GK2}.
\end{proof}

In the following let $\tilde{Q}^{f}(V,W)$ denote the square of the area of the parallelogram with sides $V$ and $W$
for $V,W\in C^{\infty}(TTM)$ given by
\begin{equation}
\tilde{Q}^{f}(V,W)=\|V\|^{2}\|W\|^{2}-\tilde{g}^{f}(V,W)^{2}.
\end{equation}

\begin{lem}\label{le:45}
Let $X,Y\in C^{\infty}(T_{p}M)$ be two orthonormal vectors in the tangent spaces $T_{p}M$ of $M$ at $p$. Then
\begin{eqnarray}
i) \  \tilde{Q}^{f}(X^{h},Y^{h})&=&f^{2},  \\
ii) \  \tilde{Q}^{f}(X^{h},Y^{v})&=&\frac{f}{\alpha}(1+g(Y,u)^{2}), \\
iii) \  \tilde{Q}^{f}(X^{v},Y^{v})&=&\frac{1}{\alpha^{2}}(1+g(Y,u)^{2}+g(X,u)^{2}).
\end{eqnarray}
\end{lem}
\begin{proof}
$i)$ The statement is a direct consequence of the definition of the Rescaled Cheeger-Gromoll Metric.

$ii)$ This is a direct consequence of
\begin{eqnarray}
 \tilde{Q}^{f}(X^{h},Y^{v})&=&\tilde{g}^{f}(X^{h},X^{h})\tilde{g}^{f}(Y^{v},Y^{v})-\tilde{g}^{f}(X^{h},Y^{v})^{2}\nonumber\\
                   &=&\frac{f}{\alpha}(1+g(Y,u)^{2}).
\end{eqnarray}

$iii)$ This last part follows from
\begin{eqnarray}
 \tilde{Q}^{f}(X^{v},Y^{v})&=&\tilde{g}^{f}(X^{v},X^{v})\tilde{g}^{f}(Y^{v},Y^{v})-\tilde{g}^{f}(X^{h},Y^{v})^{2}\nonumber\\
                   &=&\frac{1}{\alpha}(1+g(X,u)^{2})\frac{1}{\alpha}(1+g(Y,u)^{2})\nonumber\\
                      &&-[\frac{1}{\alpha^{2}}(g(X,Y)+g(X,u)g(Y,u))]^{2}        \nonumber\\
                   &=&\frac{1}{\alpha^{2}}(1+g(Y,u)^{2}+g(X,u)^{2}).
\end{eqnarray}
\end{proof}

Let $\tilde{G}^{f}$ be the $(2,0)-$tensor on the tangent bundle $TM$ given by
\begin{equation}
\tilde{G}^{f}(V,W)\mapsto \tilde{g}^{f}(\tilde{R}^{f}(V,W)W,V)
\end{equation}
for $V,W\in C^{\infty}(TTM)$.

\begin{lem}\label{le:46}
Let $X,Y\in C^{\infty}(T_{p}M)$ be two orthonormal vectors in the tangent spaces $T_{p}M$ of $M$ at $p$. Then
\begin{eqnarray}
i) \  \tilde{G}^{f}(X^{h},Y^{h})&=&\frac{1}{f} K(X,Y)-\frac{3}{4\alpha^{2}f^{2}}|R(X,Y)u|^{2}+\tilde{L}^{f}(X,Y),  \\
ii) \  \tilde{G}^{f}(X^{h},Y^{v})&=&\frac{1}{4\alpha^{2} f^{2}}|R(u,Y)X|^{2}, \\
iii) \  \tilde{G}^{f}(X^{v},Y^{v})&=&\frac{1+\alpha+\alpha^{2}}{\alpha^{2}}\tilde{Q}^{f}(X^{v},Y^{v})
                    -\frac{2+\alpha}{\alpha^{3}}(g(X,u)^{2}+g(Y,u)^{2}).\\
\end{eqnarray}
\end{lem}
\begin{proof}
$i)$ The statement follows by
\begin{eqnarray}
\alpha\tilde{G}^{f}(X^{h},Y^{h})&=&\alpha\tilde{g}^{f}(\tilde{R}^{f}(X^{h},Y^{h})Y^{h},X^{h})   \nonumber\\
                  &=& \tilde{g}^{f}\Big(\nabla_{X}(\nabla_{Y}Y+A_{f}(Y,Y))^{h},X^{h}\Big)\nonumber\\
                                  &&+\tilde{g}^{f}\Big(A_{f}(X,\nabla_{Y}Y+A_{f}(Y,Y))^{h},X^{h}\Big)\nonumber\\
                           && -\tilde{g}^{f}\Big(\nabla_{Y}(\nabla_{X}Y+A_{f}(X,Y))^{h},X^{h}\Big)\nonumber\\
                            &&-\tilde{g}^{f}\Big(A_{f}(Y,\nabla_{X}Y+A_{f}(X,Y))^{h},X^{h}\Big)\nonumber\\
                           && -\tilde{g}^{f}\Big((\nabla_{[X,Y]}Y)^{h}+A_{f}([X,Y],Y)^{h},X^{h}\Big)\nonumber\\
                             &&+\tilde{g}^{f}\Big(\frac{1}{2\alpha f(p)}\Big(R(u,R(X,Y)u)Y\Big)^{h},X^{h}\Big)\nonumber\\
                              && +\tilde{g}^{f}\Big(\frac{1}{4\alpha f(p)}\Big(R(u,R(X,Y)u)Y\Big)^{h},X^{h}\Big)\nonumber\\
                              &&-\tilde{g}^{f}\Big(\frac{1}{4\alpha f(p)}(R(u,R(Y,Y)u)X)^{h},X^{h}\Big)\nonumber\\
                          &=&\frac{1}{f} K(X,Y)-\frac{3}{4\alpha^{2}}|R(X,Y)u|^{2}+\tilde{L}^{f}(X,Y).
\end{eqnarray}

The properties of the Riemann curvature tensor give
\begin{equation}
g(R(u,R(X,Y)u)Y,X)=-|R(X,Y)u|^{2},
\end{equation}
from which the result follows.

$ii)$ The statement follows by
\begin{eqnarray}
\alpha^{2}\tilde{G}^{f}(X^{h},Y^{v})&=&\alpha^{2}\tilde{g}^{f}(\tilde{R}^{f}(X^{h},Y^{v})Y^{v},X^{h})   \nonumber\\
                   &=&-\alpha^{2}\tilde{g}^{f}\Big(-\frac{1}{2\alpha f}\Big(R(Y,Z)X\Big)^{h} ,X^{h}\Big) \nonumber\\
                    &&-\alpha^{2}\tilde{g}^{f}\Big(-\frac{1}{4\alpha^{2} f^{2}}\Big(R(u,Y)R(u,Z)X\Big)^{h} ,X^{h}\Big) \nonumber\\
                   &&+\alpha^{2}\tilde{g}^{f}\Big(\frac{1}{2\alpha^{2} f}g(Y,u)(R(u,Y)X)^{h},X^{h}\Big)\nonumber\\
                   &&-\alpha^{2}\tilde{g}^{f}\Big(\frac{1}{2\alpha^{2} f}g(Y,u)(R(u,Y)X)^{h},X^{h}\Big)\nonumber\\
                  &=&\frac{1}{4 f^{2}}|R(u,Y)X|^{2}.
\end{eqnarray}

$iii)$ In the last case we have
\begin{eqnarray}
\tilde{G}^{f}(X^{v},Y^{v})&=&\tilde{g}^{f}(\tilde{R}^{f}(X^{v},Y^{v})Y^{v},X^{v}) \nonumber\\
                 &&+\frac{\alpha+2}{\alpha^{2}}(\tilde{g}^{f}(X^{v},Y^{v})g(Y,u)g(X,u)-\tilde{g}^{f}(Y^{v},Y^{v})g(X,u)^{2})  \nonumber\\
                 &&+\frac{1+\alpha+\alpha^{2}}{\alpha^{2}}(\tilde{g}^{f}(Y^{v},Y^{v})\tilde{g}^{f}(X^{v},X^{v})
                           -\tilde{g}^{f}(X^{v},Y^{v}))   \nonumber\\
                 &&+\frac{\alpha+2}{\alpha^{2}}(g(X,u)g(Y,u)\tilde{g}^{f}(X^{v},Y^{v})
                    -g(Y,u)^{2}\tilde{g}^{f}(X^{v},X^{v}) )  \nonumber\\
                 &=&\frac{1+\alpha+\alpha^{2}}{\alpha^{2}}\tilde{Q}^{f}(X^{v},Y^{v})
                    -\frac{2+\alpha}{\alpha^{3}}(g(X,u)^{2}+g(Y,u)^{2}).
\end{eqnarray}
\end{proof}
\begin{prop}\label{pr: 47}
Let $(M, g)$ be a Riemannian manifold and $TM$ be its tangent bundle equipped with the rescaled Cheeger-Gromoll
metric $\tilde{g}^{f}$. Then the sectional curvature $\tilde{K}^{f}$ of $(TM,\tilde{g}^{f})$ satisfy the following:
\begin{eqnarray}
i) \ \tilde{K}^{f}(X^{h},Y^{h})&=&\frac{1}{f^{3}} K(X,Y)-\frac{3}{4\alpha f^{4}}|R(X,Y)u|^{2}
               +\frac{1}{f^{2}}\tilde{L}^{f}(X,Y),  \\
ii) \  \tilde{K}^{f}(X^{h},Y^{v})&=&\frac{1}{4\alpha f^{3}}\frac{|R(u,Y)X|^{2}}{(1+g(Y,u)^{2})}, \\
iii) \ \tilde{K}^{f}(X^{v},Y^{v})&=&\frac{1-\alpha}{\alpha^{2}}+\frac{2+\alpha}{\alpha}\frac{1}{(1+g(Y,u)^{2}+g(X,u)^{2})}.\\
\end{eqnarray}
\end{prop}
\begin{proof}
 The division of $\tilde{G}^{f}(X^{i},Y^{j})$ by $\tilde{Q}^{f}(X^{i},Y^{j})$ for $i,j\in\{h,v\}$ gives the result.
\end{proof}
\begin{prop}\label{pr: 48}
Let $(M, g)$ be a Riemannian manifold of constant sectional curvature $\kappa$ .Let $TM$ be its tangent bundle equipped with
the rescaled Cheeger-Gromoll metric $\tilde{g}^{f}$. Then the sectional curvature $\tilde{K}^{f}$ of $(TM,\tilde{g}^{f})$
satisfy the following:
\begin{eqnarray}
i) \ \tilde{K}^{f}(X^{h},Y^{h})&=&\frac{1}{f^{3}}\kappa -\frac{3\kappa^{2}}{4\alpha f^{4}}(g(u,X)^{2}+g(u,Y)^{2})
               +\frac{1}{f^{2}}\tilde{L}^{f}(X,Y),  \\
ii) \  \tilde{K}^{f}(X^{h},Y^{v})&=&\frac{1}{4\alpha f^{3}}\frac{\kappa^{2}g(X,u)^{2}}{(1+g(Y,u)^{2})}, \\
iii) \ \tilde{K}^{f}(X^{v},Y^{v})&=&\frac{1-\alpha}{\alpha^{2}}+\frac{2+\alpha}{\alpha}\frac{1}{(1+g(Y,u)^{2}+g(X,u)^{2})}.\\
\end{eqnarray}
for any orthonormal vectors $X,Y\in T_{p}M$.
\end{prop}
\begin{proof}
 This is a simple calculation using the special form of the curvature tensor.
\end{proof}

For a given point $(p,u)\in TM$ with $u\neq0$. Let $\{e_{1},\cdots,e_{m}\}$ be an orthonormal basis for the tangent space $T_{p}M$
of $M$ at $p$ such that $e_{1}=\frac{u}{|u|}$, where $|u|$ is the norm of $u$ with respect to the metric $g$ on $M$. Then for
 $i\in\{1,\cdots,m\}$ and $k\in\{2,\cdots,m\}$ define the horizontal and vertical lifts by $t_{i}=e_{i}^{h}$, $t_{m+1}=e_{1}^{v}$
 and $t_{m+k}=\sqrt{\alpha}e_{k}^{v}$. Then $\{t_{1},\cdots,t_{2m}\}$ is an orthonormal basis of the tangent space $T_{(p,u)}M$
 with respect to the rescaled Cheeger-Gromoll metric.
\begin{lem}\label{le:411}
Let $(p,u)$ be a point on $TM$ and $\{t_{1},\cdots,t_{2m}\}$ be an orthonormal basis of the tangent space $T_{(p,u)}M$ as above.
Then the sectional curvature $\tilde{K}^{f}$ satisfy the following equations
\begin{eqnarray}
\tilde{K}^{f}(t_{i},t_{j})&=&\frac{1}{f^{3}}K(e_{i},e_{j}) -\frac{3}{4\alpha f^{4}}|R(e_{i},e_{j})u|^{2}
               +\frac{1}{f^{2}}\tilde{L}^{f}(X,Y),  \\
\tilde{K}^{f}(t_{i},t_{m+1})&=& 0, \\
\tilde{K}^{f}(t_{i},t_{m+k})&=& \frac{1}{4f^{3}}|R(u,e_{k})e_{i}|^{2} \\
\tilde{K}^{f}(t_{m+1},t_{m+k})&=& \frac{3}{\alpha^{2}} \\
\tilde{K}^{f}(t_{m+k},t_{m+l})&=& \frac{\alpha^{2}+\alpha+1}{\alpha^{2}}  \\
\end{eqnarray}
for $i,j\in\{1,\cdots,m\}$ and $k,l\in\{2,\cdots,m\}$.
\end{lem}

\begin{prop}\label{pr: 412}
Let $(M, g)$ be a Riemannian manifold with scalar curvature $S$. Let $TM$ be its tangent bundle equipped with
the rescaled Cheeger-Gromoll metric $\tilde{g}^{f}$ and $(p,u)$ be a point on $TM$. Then the scalar curvature $\tilde{S}^{f}$
of $(TM,\tilde{g}^{f})$ satisfy the following:
Then
\begin{eqnarray}
\tilde{S}^{f}_{(p,u)}&=& S_{p}+\frac{2\alpha-3}{4\alpha f^{4}}\sum_{i,j=1}^{m}|R(e_{i},e_{j})u|^{2}+\frac{1}{f^{2}}\sum_{i,j=1}^{m}\tilde{L}^{f}(X,Y)  \nonumber\\
                                &&+\frac{m-1}{\alpha^{2}}[6+(m-2)(\alpha^{2}+\alpha+1)].
\end{eqnarray}
\end{prop}
\begin{proof}
Let $\{t_{1},\cdots,t_{2m}\}$ be an orthonormal basis of the tangent space $T_{(p,u)}TM$ as above. By the definition of the scalar
curvature we know that
\begin{eqnarray}
\tilde{S}^{f}&=&\sum_{i,j=1}^{m}\tilde{K}^{f}(t_{i},t_{j})\nonumber\\
             &=&2\sum_{i,j=1,i<j}^{m}\tilde{K}^{f}(t_{i},t_{j})+2\sum_{i,j=1}^{m}\tilde{K}^{f}(t_{i},t_{m+j})
                 +2\sum_{i,j=1,i<j}^{m}\tilde{K}^{f}(t_{m+i},t_{m+j}) \nonumber\\
             &=&\sum_{i\neq j}^{m}\tilde{K}^{f}(t_{i},t_{j})-\frac{-3}{4\alpha f^{2}}\sum_{i,j=1}^{m}\tilde{R}^{f}(t_{i},t_{j})\nonumber\\
               &&+\frac{1}{2}\sum_{i,j=1}^{m}\tilde{R}^{f}(t_{i},t_{j})+2\sum_{i=2}^{m}\frac{3}{\alpha^{2}}
                  +\sum_{i,j=1,i\neq j}^{m}\frac{\alpha^{2}+\alpha+1}{\alpha^{2}}        \nonumber\\
             &=& S+\frac{2\alpha-3}{4\alpha f^{4}}\sum_{i,j=1}^{m}|R(e_{i},e_{j})u|^{2}  \nonumber\\
                             &&+\frac{1}{f}\sum_{i,j=1}^{m}\tilde{L}^{f}(X,Y)+\frac{m-1}{\alpha^{2}}[6+(m-2)(\alpha^{2}+\alpha+1)].
\end{eqnarray}

For the fact that
\begin{equation}
\sum_{i,j=1}^{m}|R(e_{i},e_{j})u|^{2}=\sum_{i,j=1}^{m}|R(u,e_{j})e_{i}|^{2}
\end{equation}
see the proof of Proposition 3.9.
\end{proof}

\section*{ Acknowledgements}
The second author was partially supported by National Science Foundation of China under Grant No.10801027, and Fok Ying Tong
Education Foundation under Grant No.121003.

\end{document}